\documentclass[11pt]{article}
\usepackage{amsmath,amsthm,amssymb}
\usepackage[matrix,arrow]{xy}

\makeatletter
\oddsidemargin\paperwidth
\advance\oddsidemargin by-\textwidth
\advance\oddsidemargin by-1in
\evensidemargin\oddsidemargin
\headsep\footskip
\advance\headsep by-\topskip
\advance\textheight by\topskip
\topmargin\paperheight
\advance\topmargin by-\textheight
\advance\topmargin by-1in
\advance\topmargin by-\headheight
\advance\topmargin by-\headsep

\renewcommand{\section}{%
	\@startsection {section}{1}{\z@}%
	{3.5ex \@plus 1ex \@minus .2ex}%
	{2.3ex \@plus.2ex}%
	{\normalfont\large\bfseries\centering}}
\renewcommand\subsection{%
	\@startsection{subsection}{2}{\z@}%
	{3.25ex\@plus 1ex \@minus .2ex}%
	{.8ex \@plus .2ex}%
	{\normalfont\bfseries}}

\renewcommand{\@seccntformat}[1]{\S\csname the#1\endcsname.\ \,}

\renewenvironment{proof}[1][\proofname]{\par
	\pushQED{\qed}%
	\normalfont \topsep6\p@\@plus6\p@\relax
	\trivlist
	\item[%
	\scshape
	\hspace*{\parindent}
	#1\@addpunct{.}]\ignorespaces
}{%
  \popQED\endtrivlist\@endpefalse
}

\newtheoremstyle{mythm}{}{}{\itshape}{\parindent}{\bfseries}{.}{5.555pt plus 1.666pt minus 1.666pt}{}
\theoremstyle{mythm}

\renewenvironment{abstract}%
	{%
	\small
	\list{\scshape\abstractname.}%
	{%
	\labelsep 7pt
	\itemindent \labelwidth
	\advance \itemindent by \labelsep
	\listparindent 0pt
	\leftmargin 0.1\paperwidth
	\rightmargin   \leftmargin
	\parsep 0pt}
	\item\relax}%
	{\endlist}
\renewcommand{\abstractname}{Abstract}

\let\@oldenumerate=\enumerate
\def\enumerate{\@oldenumerate\def\makelabel##1{\hss ##1}}
\let\@olditemize=\itemize
\def\itemize{\@olditemize\def\makelabel##1{\hss ##1}}

\@addtoreset{equation}{section}

\def\@address{}
\def\@email{}
\def\@msc{}
\def\@mscyear{}
\newcommand{\address}[1]{\gdef\@address{#1}}
\newcommand{\email}[1]{\gdef\@email{#1}}
\newcommand{\msc}[2][2000]{\gdef\@mscyear{#1}\gdef\@msc{#2}}

\def\@org{%
\footnotetext{%
\ifx\@address\@empty\else \hspace*{-1em}{\scshape\@address}\\\fi%
\ifx\@email\@empty\else \hspace*{\footnotesep}\textit{E-mail address}: \texttt{\@email}\\\fi
\ifx\@msc\@empty\else\hspace*{\footnotesep}\@mscyear\ \textit{Mathematics Subject Classification}. \@msc.\fi}}

\let\old@maketitle=\@maketitle
\def\@maketitle{%
\old@maketitle
\ifx\@address\@empty
\ifx\@email\@empty
\ifx\@msc\@empty
\else\@org\fi
\else\@org\fi
\else\@org\fi}

\let\oldmaketitle=\maketitle
\def\maketitle{%
\let\tmp@author=\@author\let\tmp@title=\@title
\oldmaketitle
\gdef\@address{}\gdef\@email{}\gdef\@mscyear{}\gdef\@msc{}%
\let\@author=\tmp@author\let\@title=\tmp@title
}

\def\ps@mystyle{%
\def\@oddhead{\def\\{\relax}\reset@font\small\hfil\ifodd\c@page\@title\else\@author\fi\hfil}
\def\@oddfoot{\reset@font\hfil\thepage\hfil}}
\makeatother

\newtheorem{thm}{Theorem}[section]
\newtheorem{defn}[thm]{Definition}
\newtheorem{prop}[thm]{Proposition}
\newtheorem{lem}[thm]{Lemma}
\newtheorem{cor}[thm]{Corollary}

\newcommand{\R}{\mathbb{R}}
\newcommand{\C}{\mathbb{C}}
\newcommand{\Z}{\mathbb{Z}}

\DeclareMathOperator{\Gr}{Gr}
\DeclareMathOperator{\ad}{ad}
\DeclareMathOperator{\Ad}{Ad}

\DeclareMathOperator{\Ker}{Ker}
\DeclareMathOperator{\re}{Re}
\DeclareMathOperator{\Lie}{Lie}
\DeclareMathOperator{\Ind}{Ind}
\DeclareMathOperator{\Hom}{Hom}

\title{Jacquet modules of\\ principal series generated by the trivial $K$-type}
\author{Noriyuki Abe}
\address{Graduate School of Mathematical Sciences, the University of Tokyo, 3--8--1 Komaba, Meguro-ku, Tokyo 153--8914, Japan.}
\email{abenori@ms.u-tokyo.ac.jp}
\msc{22E47}
\date{}

\begin{document}
\maketitle

\begin{abstract}
We propose a new approach for the study of the Jacquet module of a Harish-Chandra module of a real semisimple Lie group.
Using this method, we investigate the structure of the Jacquet module of principal series representation generated by the trivial $K$-type.
\end{abstract}

\section{Introduction}\label{sec:Introduction}
Let $G$ be a real semisimple Lie group.
By Casselman's subrepresentation theorem, any irreducible admissible representation $U$ is realized as a subrepresentation of a certain non-unitary principal series representation.
Such an embedding is a powerful tool to study an irreducible admissible representation but the subrepresentation theorem dose not tell us how it can be realized.

Casselman~\cite{MR562655} introduced the Jacquet module $J(U)$ of $U$.
This important object retains all information of embeddings given by the subrepresentation theorem.
For example, Casselman's subrepresentation theorem is equivalent to $J(U)\ne 0$.
However the structure of $J(U)$ is very intricate and difficult to determine.

In this paper we give generators of the Jacquet module of a principal series representation generated by the trivial $K$-type.
Let $\Z$ be the ring of integers,
$\mathfrak{g}_0$ the Lie algebra of $G$,
$\theta$ a Cartan involution of $\mathfrak{g}_0$,
$\mathfrak{g}_0 = \mathfrak{k}_0\oplus\mathfrak{a}_0\oplus\mathfrak{n}_0$ the Iwasawa decomposition of $\mathfrak{g}_0$,
$\mathfrak{m}_0$ the centralizer of $\mathfrak{a}_0$ in $\mathfrak{k}_0$,
$W$ the little Weyl group for $(\mathfrak{g}_0,\mathfrak{a}_0)$,
$e\in W$ the unit element of $W$,
$\Sigma$ the restricted root system for $(\mathfrak{g}_0,\mathfrak{a}_0)$,
$\mathfrak{g}_{0,\alpha}$ the root space for $\alpha\in\Sigma$,
$\Sigma^+$ the positive system of $\Sigma$ such that $\mathfrak{n}_0 = \bigoplus_{\alpha\in\Sigma^+}\mathfrak{g}_{0,\alpha}$, 
$\rho = \sum_{\alpha\in\Sigma^+}(\dim \mathfrak{g}_{0,\alpha}/2)\alpha$,
$\mathcal{P} = \{\sum_{\alpha\in\Sigma^+}n_\alpha\alpha\mid n_\alpha\in\Z\}$,
$\mathcal{P}^+ = \{\sum_{\alpha\in\Sigma^+}n_\alpha\alpha \mid n_\alpha\in\Z,\ n_\alpha\ge 0\}$ and 
$U(\lambda)$ the principal series representation with an infinitesimal character $\lambda$ generated by the trivial $K$-type.
In this paper we prove the following theorem.
\begin{thm}[Theorem~\ref{thm:definition of boundary value map for G/K}, Theorem~\ref{thm:relations of v_i}]\label{thm:generators and relations}
Assume that $\lambda$ is regular.
Set $\mathcal{W}(w) = \{w'\in W\mid w\lambda - w'\lambda\in2\mathcal{P}^+\}$ for $w\in W$.
Then there exist generators $\{v_w\mid w\in W\}$ of $J(U(\lambda))$ such that
\[
	\begin{cases}
	\text{$(H - (\rho + w\lambda))v_w \in \sum_{w'\in \mathcal{W}(w)}U(\mathfrak{g})v_{w'}$ for all $H\in \mathfrak{a}_0$},\\
	\text{$Xv_w \in \sum_{w'\in \mathcal{W}(w)}U(\mathfrak{g})v_{w'}$ for all $X\in \mathfrak{m}_0\oplus\theta(\mathfrak{n}_0)$}.
	\end{cases}
\]
\end{thm}
We enumerate $W = \{w_1,w_2,\dots,w_r\}$ such that $\re w_1\lambda\ge \re w_2\lambda\ge \dots \ge \re w_r\lambda$.
Set $V_i = \sum_{j\ge i}U(\mathfrak{g})v_{w_j}$.
Then by Theorem~\ref{thm:generators and relations} we have the surjective map $M(w_i\lambda)\to V_i/V_{i + 1}$ where $M(w_i\lambda)$ is a generalized Verma module (See Definition~\ref{defn:generalized Verma module}).
This map is isomorphic.
Namely we can prove the following theorem.
\begin{thm}[Theorem~\ref{thm:Main theorem for regular case}]\label{thm:existence of filtration}
There exists a filtration $J(U(\lambda)) = V_1 \supset V_2\supset \dots \supset V_{r + 1} = 0$ of $J(U(\lambda))$ such that $V_i / V_{i + 1} \simeq M(w_i\lambda)$.
Moreover if $w\lambda - \lambda \not\in 2\mathcal{P}$ for $w\in W\setminus\{e\}$ then $J(U(\lambda)) \simeq \bigoplus_{w\in W} M(w\lambda)$.
\end{thm}
This theorem does not need the assumption that $\lambda$ is regular.
In the case of $G$ is split and $U(\lambda)$ is irreducible, Collingwood~\cite{MR1089304} proved Theorem~\ref{thm:existence of filtration}.

For example, we obtain the following in case of $\mathfrak{g}_0 = \mathfrak{sl}(2,\R)$:
Choose a basis $\{H,E_+,E_-\}$ of $\mathfrak{g}_0$ such that $\R H = \mathfrak{a}_0$, $\R E_+ = \mathfrak{n}_0$, $[H,E_\pm] = \pm 2E_\pm$ and $E_- = \theta(E_+)$.
Then $\Sigma^+ = \{2\alpha\}$ where $\alpha(H) = 1$.
Let $\lambda = r\alpha$ for $r\in \C$.
We may assume $\re r \ge 0$.
By Theorems~\ref{thm:generators and relations} and \ref{thm:existence of filtration}, we have the exact sequence
\[
\xymatrix{
	0 \ar[r] & M(-r\alpha) \ar[r] & J(U(r\alpha)) \ar[r] & M(r\alpha) \ar[r] & 0.
}
\]
Consider the case $\lambda$ is integral, i.e., $2r\in \Z$.
If $r\not\in \Z$ then this sequence splits by Theorem~\ref{thm:existence of filtration}.
On the other hand, if $r\in \Z$ then by the direct calculation using the method introduced in this paper we can show it does not split.
Notice that $U(r\alpha)$ is irreducible if and only if $r\in\Z$.
Then we have the following;
if $\lambda$ is integral then $J(U(\lambda))$ is isomorphic to the direct sum of generalized Verma modules if and only if $U(\lambda)$ is reducible.

Our method is based on the paper of Kashiwara and Oshima~\cite{MR0482870}.
In Section~\ref{sec:Jacquet modules and fundamental properties} we show fundamental properties of Jacquet modules and introduce a certain extension of the universal enveloping algebra.
An analog of the theory of Kashiwara and Oshima is established in Section~\ref{sec:construction of special elements}.
In Section~\ref{sec:Structure of Jacquet modules (regular case)} we prove our main theorem in the case of a regular infinitesimal character using the result of Section~\ref{sec:construction of special elements}.
We complete the proof in Section~\ref{sec:Structure of Jacquet modules (singular case)} using the translation principle.

\subsection*{Acknowledgments}
The author is grateful to his advisor Hisayosi Matumoto for his advice and support.
He would also like to thank Professor Toshio Oshima for his comments.

\subsection*{Notations}
Throughout this paper we use the following notations.
As usual we denote the ring of integers, the set of non-negative integers, the set of positive integers, the real number field and the complex number field by $\Z,\Z_{\ge 0},\Z_{> 0},\R$ and $\C$ respectively.
Let $\mathfrak{g}_0$ be a real semisimple Lie algebra.
Fix a Cartan involution $\theta$ of $\mathfrak{g}_0$.
Let $\mathfrak{g}_0 = \mathfrak{k}_0\oplus \mathfrak{s}_0$ be the decomposition of $\mathfrak{g}_0$ into the $+1$ and $-1$ eigenspaces for $\theta$.
Take a maximal abelian subspace $\mathfrak{a}_0$ of $\mathfrak{s}_0$ and let $\mathfrak{g}_0 = \mathfrak{k}_0\oplus \mathfrak{a}_0\oplus \mathfrak{n}_0$ be the corresponding Iwasawa decomposition of $\mathfrak{g}_0$.
Set $\mathfrak{m}_0 = \{X\in \mathfrak{k}_0\mid \text{$[H,X] = 0$ for all $H\in \mathfrak{a}_0$}\}$.
Then $\mathfrak{p}_0 = \mathfrak{m}_0\oplus \mathfrak{a}_0 \oplus \mathfrak{n}_0$ is a minimal parabolic subalgebra of $\mathfrak{g}_0$.
Write $\mathfrak{g}$ for the complexification of $\mathfrak{g}_0$ and $U(\mathfrak{g})$ for the universal enveloping algebra of $\mathfrak{g}$.
% and $\{U_n(\mathfrak{g})\}_{n\in\Z_{\ge 0}}$ be the standard filtration of $U(\mathfrak{g})$.
%Analogous notations are used for other Lie algebras.
We apply analogous notations to other Lie algebras.

Set $\mathfrak{a}^* = \Hom_\C(\mathfrak{a},\C)$ and $\mathfrak{a}_0^* = \Hom_\R(\mathfrak{a}_0,\R)$.
Let $\Sigma\subset\mathfrak{a}^*$ be the restricted root system for $(\mathfrak{g},\mathfrak{a})$ and $\mathfrak{g}_\alpha$ the root space for $\alpha\in\Sigma$.
Let $\Sigma^+$ be the positive root system determined by $\mathfrak{n}$, i.e., $\mathfrak{n} = \bigoplus_{\alpha\in\Sigma^+}\mathfrak{g}_\alpha$.
$\Sigma^+$ determines the set of simple roots $\Pi = \{\alpha_1,\alpha_2,\dots,\alpha_l\}$.
%We fix an order of $\mathfrak{a}_0^*$ such that $\Sigma^+ = \{\alpha\in \Sigma \mid \alpha > 0\}$.
We define the total order on $\mathfrak{a}_0^*$ by the following; for $c_i,d_i\in\R$ we define $\sum_i c_i \alpha_i > \sum_i d_i \alpha_i$ if and only if there exists an integer $k$ such that $c_1 = d_1,\dots,c_k = d_k$ and $c_{k + 1} > d_{k + 1}$.
Let $\{H_1,H_2,\dots,H_l\}$ be the dual basis of $\{\alpha_i\}$.
Write $W$ for the little Weyl group for $(\mathfrak{g}_0,\mathfrak{a}_0)$ and $e$ for the unit element of $W$.
Set 
$\mathcal{P} = \{\sum_{\alpha\in\Sigma^+}n_\alpha\alpha\mid n_\alpha\in\Z\}$, 
$\mathcal{P}^+ = \{\sum_{\alpha\in\Sigma^+}n_\alpha\alpha\mid n_\alpha\in\Z_{\ge 0}\}$
and $\mathcal{P}^{++} = \mathcal{P}^+\setminus\{0\}$.
Let $m$ be a dimension of $\mathfrak{n}$.
Fix a basis $E_1,E_2,\dots,E_m$ of $\mathfrak{n}$ such that each $E_i$ is a restricted root vector.
Let $\beta_i$ be a restricted root vector such that $E_i\in \mathfrak{g}_{\beta_i}$.
For $\mathbf{n} = (\mathbf{n}_i)\in\Z_{\ge 0}^m$ we denote $E_1^{\mathbf{n}_1}E_2^{\mathbf{n}_2}\dotsm E_m^{\mathbf{n}_m}$ by $E^{\mathbf{n}}$.

For $x = (x_1,x_2,\dots,x_n)\in\Z_{\ge 0}^n$, we write $\lvert x\rvert = x_1 + x_2 + \dots + x_n$ and $x! = x_1!\,x_2!\dotsm x_n!$.

%For sets $A$ and $\Lambda$, we denote the set of maps from $\Lambda$ to $A$ by $A^\Lambda$.
%We often write $f_\lambda$ instead of $f(\lambda)$ for $f\in A^\Lambda$ and $\lambda\in\Lambda$.

For a $\C$-algebra $R$, let $M(r,r',R)$ be the space of $r\times r'$ matrices with entries in $R$ and $M(r,R) = M(r,r,R)$.
Write $1_r\in M(r,R)$ for the identity matrix.

\section{Jacquet modules and fundamental properties}\label{sec:Jacquet modules and fundamental properties}
\begin{defn}[Jacquet module]
Let $U$ be a $U(\mathfrak{g})$-module.
Define modules $\widehat{J}(U)$ and $J(U)$ by
\begin{align*}
\widehat{J}(U) & = \varprojlim_k U/\mathfrak{n}^k U,\\
J(U) & = \widehat{J}(U)_{\text{\normalfont $\mathfrak{a}$-finite}} = \{u\in \widehat{J}(U)\mid \dim U(\mathfrak{a})u < \infty\}.
\end{align*}
We call $J(U)$ the Jacquet module of $U$.
\end{defn}

Set $\widehat{\mathcal{E}}(\mathfrak{n}) = \varprojlim_k U(\mathfrak{n})/\mathfrak{n}^kU(\mathfrak{n})$.
\begin{prop}\label{prop:algebraic property of E(n)}
\begin{enumerate}
\item The $\C$-algebra $\widehat{\mathcal{E}}(\mathfrak{n})$ is right and left Noetherian.
\item The $\C$-algebra $\widehat{\mathcal{E}}(\mathfrak{n})$ is flat over $U(\mathfrak{n})$.
\item If $U$ is a finitely generated $U(\mathfrak{n})$-module then $\varprojlim_k U/\mathfrak{n}^kU = \widehat{\mathcal{E}}(\mathfrak{n})\otimes_{U(\mathfrak{n})}U$.
\item Let $S = (S_k)$ be an element of $M(r,\mathfrak{n}\widehat{\mathcal{E}}(\mathfrak{n}))$ and $(a_n)\in \C^{\Z_{\ge 0}}$.
Define $\sum_{n = 0}^\infty a_nS^n = (\sum_{n = 0}^k a_n S_k^n)_k$.
Then $\sum_{n = 0}^\infty a_nS^n \in M(r,\widehat{\mathcal{E}}(\mathfrak{n}))$.
\end{enumerate}
\end{prop}
\begin{proof}
Since Stafford and Wallach~\cite[Theorem~2.1]{MR656493} show that $\mathfrak{n}U(\mathfrak{n})\subset U(\mathfrak{n})$ satisfies the Artin-Rees property, the usual argument of the proof for commutative rings can be applicable to prove (1), (2) and (3).
(4) is obvious.
\end{proof}

\begin{cor}
Let $S$ be an element of $M(r,\widehat{\mathcal{E}}(\mathfrak{n}))$ such that $S - 1_r \in M(r,\mathfrak{n}\widehat{\mathcal{E}}(\mathfrak{n}))$.
Then $S$ is invertible.
\end{cor}
\begin{proof}
Set $T = 1_r - S$.
By Proposition~\ref{prop:algebraic property of E(n)}, $R = \sum_{n = 0}^\infty T^n\in M(r,\widehat{\mathcal{E}}(\mathfrak{n}))$.
Then $SR = RS = 1_r$.
\end{proof}

We can prove the following proposition in a similar way to that of Goodman and Wallach~\cite[Lemma~2.2]{MR597811}.
For the sake of completeness we give a proof.

\begin{prop}\label{prop:generalized result of Goodman and Wallach}
Let $U$ be a $U(\mathfrak{a}\oplus\mathfrak{n})$-module which is finitely generated as a $U(\mathfrak{n})$-module.
Assume that every element of $U$ is $\mathfrak{a}$-finite.
For $\mu\in\mathfrak{a}^*$ set
\[
	U_\mu = \{u\in U\mid \text{For all $H\in\mathfrak{a}$ there exists a positive integer $N$ such that $(H - \mu(H))^Nu = 0$}\}.
\]
Then
\[
	\widehat{J}(U) \simeq \prod_{\mu\in\mathfrak{a}^*}U_\mu.
\]
\end{prop}
\begin{proof}
For $k\in\Z_{>0}$ put $S_k = \{\mu\in\mathfrak{a}^* \mid U_\mu\ne 0,\ U_\mu\not\subset \mathfrak{n}^kU\}$.
Since $U$ is finitely generated, $\dim U/\mathfrak{n}^kU < \infty$.
Therefore $S_k$ is a finite set.
Define a map $\varphi\colon \prod_{\mu\in\mathfrak{a}^*} U_\mu\to \widehat{J}(U)$ by
\[
	\varphi((x_\mu)_{\mu\in\mathfrak{a}^*}) = \left(\sum_{\mu\in S_k} x_\mu \pmod{\mathfrak{n}^kU}\right)_k.
\]

First we show that $\varphi$ is injective.
Assume $\varphi((x_\mu)_{\mu\in\mathfrak{a}^*}) = 0$.
We have $\sum_{\mu\in S_k}x_\mu \in \mathfrak{n}^kU$ for all $k\in\Z_{> 0}$.
Since $\mathfrak{n}^kU$ is $\mathfrak{a}$-stable and $S_k$ is a finite set, $x_\mu\in\mathfrak{n}^kU$ for all $\mu\in\mathfrak{a}^*$, thus we have $x_\mu = 0$.%  $U$ is finitely generated $U(\mathfrak{n})$-module.

We have to show that $\varphi$ is surjective.
Let $x = (x_k\pmod{\mathfrak{n}^k U})_k$ be an element of $\widehat{J}(U)$.
Since every element of $U$ is $\mathfrak{a}$-finite, we have $U = \bigoplus_{\mu\in\mathfrak{a}^*}U_\mu$.
Let $p_\mu\colon U\to U_\mu$ be the projection.
The $U(\mathfrak{n})$-module $U$ is finitely generated and therefore for all $\mu\in\mathfrak{a}^*$ there exists a positive integer $k_\mu$ such that $\mathfrak{n}^{k_\mu}U\cap U_\mu = 0$.
Notice that if $i,i' > k_\mu$ then $p_\mu(x_i) = p_\mu(x_{i'})$.
Hence we have $\varphi((p_\mu(x_{k_\mu}))_{\mu\in\mathfrak{a}^*}) = x$.
\end{proof}

We define an $(\mathfrak{a}\oplus\mathfrak{n})$-representation structure of $U(\mathfrak{n})$ by $(H + X)(u) = Hu - uH + Xu$ for $H\in \mathfrak{a}$, $X\in \mathfrak{n}$, $u\in U(\mathfrak{n})$.
Then $U(\mathfrak{n})$ is a $U(\mathfrak{a}\oplus\mathfrak{n})$-module.
By Proposition~\ref{prop:generalized result of Goodman and Wallach} $\widehat{\mathcal{E}}(\mathfrak{n}) = \prod_{\mu\in\mathfrak{a}^*}U(\mathfrak{n})_\mu$.
The following results are corollaries of Proposition~\ref{prop:generalized result of Goodman and Wallach}.

\begin{cor}\label{cor:E(n) as vector space}
A linear map
\[
	\begin{array}{ccc}
	\C[[X_1,X_2,\dots,X_m]] & \longrightarrow & \widehat{\mathcal{E}}(\mathfrak{n})\\
	\sum_{\mathbf{n}\in\Z_{\ge 0}^m} a_{\mathbf{n}} X^{\mathbf{n}}
	& \longmapsto &
	\left(\sum_{\lvert\mathbf{n}\rvert \le k}a_{\mathbf{n}} E^{\mathbf{n}}\pmod{\mathfrak{n}^kU(\mathfrak{n})}\right)_k
	\end{array}
\]
is bijective, where $X^{\mathbf{n}} = X_1^{\mathbf{n}_1}X_2^{\mathbf{n}_2}\dotsm X_m^{\mathbf{n}_m}$ for $\mathbf{n} = (\mathbf{n}_1,\mathbf{n}_2,\dots,\mathbf{n}_m)\in\Z_{\ge 0}^m$.
\end{cor}
\begin{proof}
By the Poincar\'e-Birkhoff-Witt theorem $\{E^{\mathbf{n}}\mid \sum_i\mathbf{n}_i\beta_i = \mu\}$ is a basis of $U(\mathfrak{n})_\mu$.
This implies the corollary since $\widehat{\mathcal{E}}(\mathfrak{n})= \prod_{\mu\in\mathfrak{a}^*}U(\mathfrak{n})_\mu$.
\end{proof}

We denote the image of $\sum_{\mathbf{n}\in\Z_{\ge 0}^m} a_{\mathbf{n}}X^{\mathbf{n}}$ under the map in Corollary~\ref{cor:E(n) as vector space} by $\sum_{\mathbf{n}\in\Z_{\ge 0}^m} a_{\mathbf{n}}E^{\mathbf{n}}$.

\begin{cor}\label{cor:Jacquet module of trivial case}
Let $U$ be a $U(\mathfrak{g})$-module which is finitely generated as a $U(\mathfrak{n})$-module.
Assume that all elements are $\mathfrak{a}$-finite.
Then $J(U) = U$.
\end{cor}
\begin{proof}
This follows from the following equation.
\[
	J(U) = \widehat{J}(U)_{\text{$\mathfrak{a}$-finite}} = \left(\prod_{\mu\in\mathfrak{a}^*}U_\mu\right)_{\text{$\mathfrak{a}$-finite}} = \bigoplus_{\mu\in\mathfrak{a}^*}U_\mu = U.
\]
\end{proof}

Put $\widehat{\mathcal{E}}(\mathfrak{g},\mathfrak{n}) = \widehat{\mathcal{E}}(\mathfrak{n})\otimes_{U(\mathfrak{n})}U(\mathfrak{g})$.
We can define a $\C$-algebra structure of $\widehat{\mathcal{E}}(\mathfrak{g},\mathfrak{n})$ by
\begin{align*}
(f\otimes 1)(1\otimes u) & = f\otimes u,\\
(1\otimes u)(1\otimes u') & = 1\otimes (uu'),\\
(f\otimes 1)(f'\otimes 1) & = (ff')\otimes 1,\\
(1\otimes X)(f\otimes 1) & = \sum_{{\mathbf{n}}\in\Z_{\ge 0}^r}\frac{1}{{\mathbf{n}}!}\frac{\partial^{\lvert{\mathbf{n}}\rvert}}{\partial E^{\mathbf{n}}}f\otimes ((\ad(E))^{\mathbf{n}})'(X),
\end{align*}
where $u,u'\in U(\mathfrak{g})$, $X\in\mathfrak{g}$, $f,f'\in\widehat{\mathcal{E}}(\mathfrak{n})$, $((\ad(E))^{\mathbf{n}})' = (-\ad(E_m))^{{\mathbf{n}}_m}
%(-\ad(E_{m - 1}))^{{\mathbf{n}}_{m - 1}}
\dotsm(-\ad(E_1))^{{\mathbf{n}}_1}$ and 
\[
	\frac{\partial^{\lvert{\mathbf{n}}\rvert}}{\partial E^{\mathbf{n}}}\left(\sum_{\mathbf{m}\in\Z_{\ge 0}^m} a_{\mathbf{m}} E^{\mathbf{m}}\right) = \sum_{\mathbf{m}\in\Z_{\ge 0}^m} a_{\mathbf{m}}\frac{{\mathbf{m}}!}{({\mathbf{m}}-{\mathbf{n}})!}E^{{\mathbf{m}}-{\mathbf{n}}}.
\]
Notice that $\widehat{\mathcal{E}}(\mathfrak{g},\mathfrak{n})\otimes_{U(\mathfrak{g})}U \simeq \widehat{\mathcal{E}}(\mathfrak{n})\otimes_{U(\mathfrak{n})}U$ as an $\widehat{\mathcal{E}}(\mathfrak{n})$-module for a $U(\mathfrak{g})$-module $U$.
By Proposition~\ref{prop:algebraic property of E(n)}, $\widehat{\mathcal{E}}(\mathfrak{g},\mathfrak{n})$ is flat over $U(\mathfrak{g})$.
Notice that if $\mathfrak{b}$ is a subalgebra of $\mathfrak{g}$ which contains $\mathfrak{n}$ then $\widehat{\mathcal{E}}(\mathfrak{n})\otimes_{U(\mathfrak{n})}U(\mathfrak{b})$ is a subalgebra of $\widehat{\mathcal{E}}(\mathfrak{g},\mathfrak{n})$.
Put $\widehat{\mathcal{E}}(\mathfrak{b},\mathfrak{n}) = \widehat{\mathcal{E}}(\mathfrak{n})\otimes_{U(\mathfrak{n})}U(\mathfrak{b})$.

Let $U$ be a $U(\mathfrak{a}\oplus\mathfrak{n})$-module such that $U = \bigoplus_{\mu\in\mathfrak{a}^*}U_\mu$. %where
%\begin{multline*}
%U_\mu = \{u\in U\mid \\
%\text{for each $H\in\mathfrak{a}$ there exists a positive integer $N$ such that $(H - \mu(H))^Nu = 0$}\}.
%\end{multline*}
Set 
\[
	V = \left\{(u_\mu)_\mu\in\prod_{\mu\in\mathfrak{a}^*}U_\mu\Biggm| \text{there exists an element $\nu \in \mathfrak{a}_0^*$ such that $u_\mu = 0$ for $\re\mu < \nu$}\right\}.
\]
Then we can define an $\mathfrak{a}$-module homomorphism
\[
	\varphi\colon \widehat{\mathcal{E}}(\mathfrak{a}\oplus\mathfrak{n},\mathfrak{n})\otimes_{U(\mathfrak{a}\oplus\mathfrak{n})}U\simeq \left(\prod_{\mu\in\mathfrak{a}^*} U(\mathfrak{n})_\mu\right)\otimes_{U(\mathfrak{n})} \left(\bigoplus_{\mu'\in\mathfrak{a}^*}U_{\mu'}\right) \to V
\]
by
$\varphi((f_\mu)_{\mu\in\mathfrak{a}^*}\otimes (u_{\mu'})_{\mu'\in\mathfrak{a}^*}) = (\sum_{\mu + \mu' = \lambda}f_\mu u_{\mu'})_{\lambda\in\mathfrak{a}^*}$.
Notice that the composition $U\to \widehat{U}\to V$ is equal to the inclusion map $U\hookrightarrow V$.

We consider the case $U = U(\mathfrak{g})$.
Define an $(\mathfrak{a}\oplus \mathfrak{n})$-module structure of $U(\mathfrak{g})$ by $(H + X)(u) = Hu - uH + Xu$ for $H\in \mathfrak{a}$, $X\in \mathfrak{n}$, $u\in U(\mathfrak{g})$.
We have a map
\begin{multline*}
\varphi\colon \widehat{\mathcal{E}}(\mathfrak{g},\mathfrak{n})\to \\
\left\{(P_\mu)_{\mu\in\mathfrak{a}^*}\in\prod_{\mu\in\mathfrak{a}^*}U(\mathfrak{g})_\mu\,\Bigg|\,\text{there exists an element $\nu \in \mathfrak{a}_0^*$ such that $P_\mu = 0$ for $\re\mu < \nu$}\right\}.
\end{multline*}
We write $\varphi(P) = (P^{(\mu)})_{\mu\in\mathfrak{a}^*}$.
Put $P^{(H,z)} = \sum_{\mu(H) = z}P^{(\mu)}$ for $z\in \C$ and $H\in \mathfrak{a}$ such that $\re\alpha(H) > 0$ for all $\alpha\in\Sigma^+$.

By the definition we have the following proposition.
\begin{prop}\label{prop:fundamental property of ^(lambda)}
\begin{enumerate}
\item Assume that $U$ is finitely generated as a $U(\mathfrak{n})$-module.
Let $\varphi\colon \widehat{\mathcal{E}}(\mathfrak{a}\oplus\mathfrak{n},\mathfrak{n})\otimes_{U(\mathfrak{n})}U\to \prod_{\mu\in\mathfrak{a}^*}U_\mu$ be an $\mathfrak{a}$-module homomorphism defined as above.
Then $\varphi$ is coincide with the map given in Proposition~\ref{prop:generalized result of Goodman and Wallach}.
In particular $\varphi$ is isomorphic.
\item We have $(PQ)^{(\lambda)} = \sum_{\mu + \mu' = \lambda}P^{(\mu)}Q^{(\mu')}$ for $P,Q\in\widehat{\mathcal{E}}(\mathfrak{g},\mathfrak{n})$ and $\lambda\in\mathfrak{a}^*$.
\item We have
\[
	\left(\sum_{\mathbf{n}\in\Z_{\ge 0}^m}a_{\mathbf{n}}E^{\mathbf{n}}\right)^{(\lambda)} = \sum_{\sum_i \mathbf{n}_i\beta_i = \lambda}a_{\mathbf{n}}E^{\mathbf{n}}
\]
for $\lambda\in\mathfrak{a}^*$.
\end{enumerate}
\end{prop}

\begin{prop}\label{prop:induced equation}
Let $U$ be a $U(\mathfrak{g})$-module which is finitely generated as a $U(\mathfrak{n})$-module.
We take generators $v_1,v_2,\dots,v_n$ of an $\widehat{\mathcal{E}}(\mathfrak{g},\mathfrak{n})$-module $\widehat{\mathcal{E}}(\mathfrak{g},\mathfrak{n})\otimes_{U(\mathfrak{g})}U$ and set $V = \sum_i U(\mathfrak{g})v_i\subset \widehat{\mathcal{E}}(\mathfrak{g},\mathfrak{n})\otimes_{U(\mathfrak{g})}U$.
Define the surjective map $\psi\colon \widehat{\mathcal{E}}(\mathfrak{g},\mathfrak{n})\otimes_{U(\mathfrak{g})}V \to \widehat{\mathcal{E}}(\mathfrak{g},\mathfrak{n})\otimes_{U(\mathfrak{g})}U$ by $\psi(f\otimes v) = fv$.
Assume that there exist weights $\lambda_i\in\mathfrak{a}^*$ and a positive integer $N$ such that $(H - \lambda_i(H))^Nv_i = 0$ for all $H\in \mathfrak{a}$ and $1 \le i \le n$.
Let $\varphi\colon \widehat{\mathcal{E}}(\mathfrak{g},\mathfrak{n})\otimes_{U(\mathfrak{g})}V \to \prod_{\mu\in\mathfrak{a}^*} V_\mu$ be the map defined as above.
Then there exists a unique map $\widehat{\mathcal{E}}(\mathfrak{g},\mathfrak{n})\otimes_{U(\mathfrak{g})}U\to \prod_{\mu\in\mathfrak{a}^*} V_\mu$ such that the diagram
\[
	\xymatrix{
	\widehat{\mathcal{E}}(\mathfrak{g},\mathfrak{n})\otimes_{U(\mathfrak{g})}V \ar[r]^(.57){\varphi} \ar[d]^\psi &
	\prod_{\mu\in\mathfrak{a}^*}V_\mu\\
	\widehat{\mathcal{E}}(\mathfrak{g},\mathfrak{n})\otimes_{U(\mathfrak{g})}U \ar@{.>}[ru]}
\]
is commutative.
\end{prop}
\begin{proof}
Set $\widehat{U} = \widehat{\mathcal{E}}(\mathfrak{g},\mathfrak{n})\otimes_{U(\mathfrak{g})}U$ and $\widehat{V} = \widehat{\mathcal{E}}(\mathfrak{g},\mathfrak{n})\otimes_{U(\mathfrak{g})}V$.
Take $f^{(i)}\in \widehat{\mathcal{E}}(\mathfrak{g},\mathfrak{n})$ and $v^{(i)}\in V$ such that $\psi(\sum_i f^{(i)}\otimes v^{(i)}) = 0$.
We have to show $\varphi(\sum_i f^{(i)}\otimes v^{(i)}) = 0$.
Since $\widehat{V} = \widehat{\mathcal{E}}(\mathfrak{n})\otimes_{U(\mathfrak{n})}V$, we may assume $f_i\in \widehat{\mathcal{E}}(\mathfrak{n})$.
We can write $f^{(i)} = (f_\mu^{(i)})_{\mu\in\mathfrak{a}^*}$ by the isomorphism $\widehat{\mathcal{E}}(\mathfrak{n}) \simeq \prod_{\mu\in\mathfrak{a}^*}U(\mathfrak{n})_\mu$.
Since $V = \bigoplus_{\mu'\in\mathfrak{a}^*}V_{\mu'}$, we can write $v_i = \sum_{\mu'\in\mathfrak{a}^*}v_{\mu'}^{(i)}$, $v_{\mu'}^{(i)}\in V_{\mu'}$.
We have to show $\sum_i\sum_{\mu + \mu' = \lambda}f^{(i)}_\mu v^{(i)}_{\mu'} = 0$ for all $\lambda\in\mathfrak{a}^*$.
Since $U$ is a finitely generated $U(\mathfrak{n})$-module we have $\widehat{U} = \varprojlim_k U/\mathfrak{n}^kU = \varprojlim_k \widehat{U}/\mathfrak{n}^k\widehat{U}$.
It is sufficient to prove $\sum_i\sum_{\mu + \mu' = \lambda}f^{(i)}_\mu v^{(i)}_{\mu'}\in \mathfrak{n}^k\widehat{U}$ for all $k\in\Z_{>0}$.

Fix $\lambda\in\mathfrak{a}^*$ and $k\in\Z_{>0}$.
We can choose an element $\nu\in\mathfrak{a}^*_0$ such that
$\bigoplus_{\re\mu \ge \nu}U(\mathfrak{n})_\mu\subset \mathfrak{n}^kU(\mathfrak{n})$.
Then $0 = \varphi(\sum_i f^{(i)}\otimes v^{(i)}) \equiv \sum_i\sum_{\re\mu < \nu}f^{(i)}_\mu v^{(i)}_{\mu'}\pmod{\mathfrak{n}^k\widehat{U}}$.
Notice that following two sets are finite.
\begin{gather*}
\{\mu \mid \text{$\re(\mu) < \nu$ and there exists an integer $i$ such that $f_\mu^{(i)} \ne 0$}\},\\
\{\mu'\mid \text{there exists an integer $i$ such that $v_{\mu'}^{(i)} \ne 0$}\}.
\end{gather*}
This implies $\sum_i\sum_{\mu + \mu' = \lambda}f_\mu^{(i)}v_{\mu'}^{(i)} \in \mathfrak{n}^k\widehat{U}$.
%Notice that $\widehat{\mathcal{E}}(\mathfrak{n}) = \prod_{\mu\in\mathfrak{a}^*}U(\mathfrak{n})_\mu$ and $V = \bigoplus_{\mu'\in\mathfrak{a}^*}V_{\mu'}$.
%Suppose $\sum_{i} (f^{(i)}_\mu)_{\mu\in\mathfrak{a}^*}\otimes (v^{(i)}_{\mu'})_{\mu'\in\mathfrak{a}^*}\in \widehat{\mathcal{E}}(\mathfrak{g},\mathfrak{n})\otimes_{U(\mathfrak{g})}V$ satisfies $\sum_i(f^{(i)}_\mu)_{\mu\in\mathfrak{a}^*}(v^{(i)}_{\mu'})_{\mu'\in\mathfrak{a}^*} = 0$ in $\widehat{\mathcal{E}}(\mathfrak{g},\mathfrak{n})\otimes_{U(\mathfrak{g})}U$.
%We have to show $\sum_{i}\sum_{\mu + \mu' = \lambda}f^{(i)}_\mu v^{(i)}_{\mu'} = 0$ for each $\lambda\in\mathfrak{a}^*$.
%
%Set $\widehat{U} = \widehat{\mathcal{E}}(\mathfrak{g},\mathfrak{n})\otimes_{U(\mathfrak{g})}U$.\
%Since $U$ is a finitely generated $U(\mathfrak{n})$-module, $\widehat{U} = \varprojlim_k U/\mathfrak{n}^kU$.
%We will show $\sum_{i}\sum_{\mu + \mu' = \lambda}f^{(i)}_\mu v^{(i)}_{\mu'} \equiv 0\pmod{\mathfrak{n}^k\widehat{U}}$ for all $k\in\Z_{>0}$.
%
%Fix $\lambda\in\mathfrak{a}^*$ and $k\in\Z_{>0}$.
%We can choose an element $r \in\mathfrak{a}_0^*$ such that
%\[
%	\bigoplus_{\re\mu \ge r}U(\mathfrak{n})_\mu \subset \mathfrak{n}^kU(\mathfrak{n}).
%\]
%Then
%\[
%	\sum_{i}(f^{(i)}_\mu)(v^{(i)}_{\mu'}) \equiv
%	\sum_{i}\sum_{\re(\mu) < r}f^{(i)}_\mu v^{(i)}_{\mu'}
%	\pmod{\mathfrak{n}^k\widehat{U}}.
%\]
%Since $\{\mu + \mu'\mid \text{$\re(\mu) < r$, $\{i\mid f_{\mu}^{(i)}v_{\mu'}^{(i)} \ne 0\} \ne\emptyset$}\}$ is a finite set, we have
%\[
%	\sum_{i}\sum_{\mu + \mu' = \lambda}f^{(i)}_\mu v^{(i)}_{\mu'} \equiv 0\pmod{\mathfrak{n}^k\widehat{U}}.
%\]
\end{proof}

The following result is a corollary of Proposition~\ref{prop:induced equation}.

\begin{cor}\label{cor:induced equation}
In the setting of Proposition~\ref{prop:induced equation}, we have the following.
Let $P_i$ $(1 \le i \le n)$ be elements of $\widehat{\mathcal{E}}(\mathfrak{g},\mathfrak{n})$ such that $\sum_{i = 1}^nP_iv_i = 0$.
Then $\sum_iP_i^{(\lambda - \lambda_i)}v_i = 0$ for all $\lambda\in\mathfrak{a}^*$.
\end{cor}

%\begin{cor}\label{cor:induced equation2}
%Let $H$ be an element of $\mathfrak{a}$ such that $\re\alpha(H) > 0$ for all $\alpha\in \Sigma^+$.
%Then $\sum_iP_i^{(H,z - \lambda_i(H))}v_i = 0$ for all $z\in \C$.
%\end{cor}

\section{Construction of special elements}\label{sec:construction of special elements}
Let $\Lambda$ be a subset of $\mathcal{P}$.
Put $\Lambda^+ = \Lambda\cap \mathcal{P}^+$ and $\Lambda^{++} = \Lambda\cap \mathcal{P}^{++}$.
We define vector spaces $U(\mathfrak{g})_\Lambda, U(\mathfrak{n})_\Lambda, \widehat{\mathcal{E}}(\mathfrak{n})_\Lambda$ and $\widehat{\mathcal{E}}(\mathfrak{g},\mathfrak{n})_\Lambda$ by
\begin{align*}
U(\mathfrak{g})_\Lambda & = \{P\in U(\mathfrak{g})\mid \text{$P^{(\mu)} = 0$ for all $\mu\not\in\Lambda$}\},\\
U(\mathfrak{n})_\Lambda & = \{P\in U(\mathfrak{n})\mid \text{$P^{(\mu)} = 0$ for all $\mu\not\in\Lambda$}\},\\
\widehat{\mathcal{E}}(\mathfrak{n})_\Lambda & = \{P\in \widehat{\mathcal{E}}(\mathfrak{n})\mid \text{$P^{(\mu)} = 0$ for all $\mu\not\in\Lambda$}\},\\
\widehat{\mathcal{E}}(\mathfrak{g},\mathfrak{n})_\Lambda & = \{P \in \widehat{\mathcal{E}}(\mathfrak{g},\mathfrak{n})\mid \text{$P^{(\mu)} = 0$ for all $\mu\not\in\Lambda$}\}.
\end{align*}
Put $(\mathfrak{n}U(\mathfrak{n}))_\Lambda = \mathfrak{n}U(\mathfrak{n})\cap U(\mathfrak{n})_\Lambda$ and $(\mathfrak{n}\widehat{\mathcal{E}}(\mathfrak{n}))_\Lambda = \mathfrak{n}\widehat{\mathcal{E}}(\mathfrak{n})\cap \widehat{\mathcal{E}}(\mathfrak{n})_\Lambda$.

We assume that $\Lambda$ is a subgroup of $\mathfrak{a}^*$.
Then $U(\mathfrak{g})_\Lambda, U(\mathfrak{n})_\Lambda, \widehat{\mathcal{E}}(\mathfrak{n})_\Lambda$ and $\widehat{\mathcal{E}}(\mathfrak{g},\mathfrak{n})_\Lambda$ are $\C$-algebras.
Let $U$ be a $U(\mathfrak{g})_\Lambda$-module which is finitely generated as a $U(\mathfrak{n})_\Lambda$-module.
Let $u_1,u_2,\dots,u_N$ be generators of $U$ as a $U(\mathfrak{n})_\Lambda$-module.
Put $u = {}^t(u_1,u_2,\dots,u_N)$, $\overline{U} = U/(\mathfrak{n}U(\mathfrak{n}))_\Lambda U$ and $\overline{u} = u\pmod{(\mathfrak{n}U(\mathfrak{n}))_\Lambda}$.
The module $\overline{U}$ has an $\mathfrak{a}$-module structure induced from that of $U$.
By the assumption we have $\dim \overline{U} < \infty$.
Let $\lambda_1,\lambda_2,\dots,\lambda_r \in \mathfrak{a}^*$ $(\re\lambda_1 \ge \re \lambda_2\ge \dots \ge \re\lambda_r)$ be eigenvalues of $\mathfrak{a}$ on $\overline{U}$ with multiplicities.
We can choose a basis $\overline{v_1},\overline{v_2},\dots,\overline{v_r}$ of $\overline{U}$ and a linear map $\overline{Q}\colon\mathfrak{a}\to M(r,\C)$ such that
\[
	\begin{cases}
	\text{$H\overline{v} = \overline{Q}(H)\overline{v}$ for all $H\in \mathfrak{a}$},\\
	\text{$\overline{Q}(H)_{ii} = \lambda_i(H)$ for all $H\in \mathfrak{a}$}, \\
	\text{if $i > j$ then $\overline{Q}(H)_{ij} = 0$ for all $H\in \mathfrak{a}$}, \\
	\text{if $\lambda_i \ne \lambda_j$ then $\overline{Q}(H)_{ij} = 0$ for all $H\in \mathfrak{a}$},
	\end{cases}
\]
where $\overline{v} = {}^t(\overline{v_1},\overline{v_2},\dots,\overline{v_r})$.
Take $\overline{A}\in M(N,r,\C)$ and $\overline{B}\in M(r,N,\C)$ such that $\overline{v} = \overline{B}\overline{u}$ and $\overline{u} = \overline{A}\overline{v}$.

Set $\widehat{U} = \widehat{\mathcal{E}}(\mathfrak{g},\mathfrak{n})_\Lambda\otimes_{U(\mathfrak{g})_\Lambda}U$.

\begin{thm}\label{thm:definition of boundary value map}
There exist matrices $A\in M(N,r,\widehat{\mathcal{E}}(\mathfrak{n})_\Lambda)$ and $B\in M(r,N,\widehat{\mathcal{E}}(\mathfrak{n})_\Lambda)$ such that the following conditions hold:
\begin{itemize}
\item There exists a linear map $Q\colon \mathfrak{a}\to M(r,U(\mathfrak{n})_\Lambda)$ such that
\[
	\begin{cases}
	\text{$Hv = Q(H)v$ for all $H\in \mathfrak{a}$},\\
	\text{$Q(H) - \overline{Q}(H)\in M(r,(\mathfrak{n}U(\mathfrak{n}))_\Lambda)$ for all $H\in\mathfrak{a}$},\\
	\text{if $\lambda_i - \lambda_j \not\in\Lambda^+$ then $Q(H)_{ij} = 0$ for all $H\in \mathfrak{a}$},\\
	\text{if $\lambda_i - \lambda_j \in \Lambda^+$ then $[H',Q(H)_{ij}] = (\lambda_i - \lambda_j)(H')Q(H)_{ij}$ for all $H,H'\in \mathfrak{a}$},
	\end{cases}
\]
where $v = Bu$.
\item We have $u = ABu$.
\item We have $A - \overline{A} \in M(N,r,(\mathfrak{n}\widehat{\mathcal{E}}(\mathfrak{n}))_\Lambda)$ and $B - \overline{B}\in M(r,N,(\mathfrak{n}\widehat{\mathcal{E}}(\mathfrak{n}))_\Lambda)$.
\end{itemize}
\end{thm}

For the proof we need some lemmas.
Put $w = \overline{B}u\in \widehat{U}^r$.

\begin{lem}\label{lem:lemma for proof of boudary value map}
For $H\in\mathfrak{a}$ there exists a matrix $R\in M(r,(\mathfrak{n}\widehat{\mathcal{E}}(\mathfrak{n}))_\Lambda)$ such that $Hw = (\overline{Q}(H) + R)w$ in $\widehat{U}^r$.
\end{lem}
\begin{proof}
Since $w \pmod{((\mathfrak{n}U(\mathfrak{n}))_\Lambda U)^r} = \overline{v}$, we have $Hw - \overline{Q}(H)w \in ((\mathfrak{n}U(\mathfrak{n}))_\Lambda U)^r$.
Hence there exists a matrix $R_1 \in M(N,r,(\mathfrak{n}U(\mathfrak{n}))_\Lambda)$ such that $Hw - \overline{Q}(H)w = R_1u$.
Similarly we can choose a matrix $S\in M(N,(\mathfrak{n}U(\mathfrak{n}))_\Lambda)$ which satisfies $u = \overline{A}w + Su$.
Put $R = R_1(1 - S)^{-1}\overline{A}$.
Then $(H - \overline{Q}(H) -  R)w = R_1u - R_1(1 - S)^{-1}\overline{A}w = 0$.
\end{proof}

\begin{lem}\label{lem:Lemma of linear algebra}
Let $\lambda\in\C$ and $Q_0,R\in M(r,\C)$.
Assume that $Q_0$ is an upper triangular matrix.
Then there exist matrices $L,T\in M(r,\C)$ such that 
\[
	\begin{cases}
	\lambda L - [Q_0,L] = T + R,\\
	\text{if $(Q_0)_{ii} - (Q_0)_{jj} \ne \lambda$ then $T_{ij} = 0$}.
	\end{cases}
\]
\end{lem}
\begin{proof}
Since $(Q_0)_{ij} = 0$ for $i > j$, we have
\begin{align*}
(\lambda L - [Q_0,L])_{ij} & = \lambda L_{ij} - \sum_{k = 1}^r ((Q_0)_{ik}L_{kj} - L_{ik}(Q_0)_{kj})\\
& = \lambda L_{ij} - \sum_{k = i}^r (Q_0)_{ik}L_{kj} + \sum_{k = 1}^jL_{ik}(Q_0)_{kj}\\
& = (\lambda - ((Q_0)_{ii} - (Q_0)_{jj}))L_{ij} - \sum_{k = i + 1}^r (Q_0)_{ik}L_{kj} + \sum_{k = 1}^{j - 1}L_{ik}(Q_0)_{kj}.
\end{align*}
Hence we can choose $L_{ij}$ and $T_{ij}$ inductively on $(j - i)$.
\end{proof}

\begin{lem}\label{lem:Main lemma of def of boundary value map}
Let $H$ be an element of $\mathfrak{a}$ such that $\alpha(H) > 0$ for all $\alpha\in \Sigma^+$.
Let $Q_0\in M(r,\C)$, $R\in M(r,(\mathfrak{n}\widehat{\mathcal{E}}(\mathfrak{n}))_\Lambda)$.
Assume $(Q_0)_{ij} = 0$ for $i > j$.
Set $\mathcal{L}^{++} = \{\lambda(H)\mid \lambda\in\Lambda^{++}\}$.
Then there exist matrices $L\in M(r,\widehat{\mathcal{E}}(\mathfrak{n})_\Lambda)$ and $T\in M(r,(\mathfrak{n}U(\mathfrak{n}))_\Lambda)$ such that
\[
	\begin{cases}
	L \equiv 1_r \pmod{(\mathfrak{n}U(\mathfrak{n}))_\Lambda},\\
	(H1_r - Q_0 - T)L = L(H1_r - Q_0 - R),\\
%	\text{$T_{ij} = 0$ \quad if $i > j$},\\
	\text{if $(Q_0)_{ii} - (Q_0)_{jj} \not\in\mathcal{L}^{++}$ then $T_{ij} = 0$},\\
	\text{if $(Q_0)_{ii} - (Q_0)_{jj}\in\mathcal{L}^{++}$ then $[H,T_{ij}] = ((Q_0)_{ii} - (Q_0)_{jj})T_{ij}$}.
	\end{cases}
\]
\end{lem}
\begin{proof}
Set $\mathcal{L} = \{\lambda(H)\mid \lambda\in\Lambda\}$ and $\mathcal{L}^+ = \{\lambda(H)\mid \lambda\in\Lambda^+\}$.
Put $f(\mathbf{n}) = \sum_i \mathbf{n}_i\beta_i$ for $\mathbf{n} = (\mathbf{n}_i) \in \Z^m$.
Set $\widetilde{\Lambda} = \{\mathbf{n}\in \Z_{\ge 0}^m\mid f(\mathbf{n}) \in \Lambda\}$.
We define the order on $\widetilde{\Lambda}$ by $\mathbf{n} < \mathbf{n}'$ if and only if $f(\mathbf{n}) < f(\mathbf{n}')$.

By Corollary~\ref{cor:E(n) as vector space}, we can write $R = \sum_{\mathbf{n}\in \widetilde{\Lambda}}R_{\mathbf{n}}E^{\mathbf{n}}$ where $R_{\mathbf{n}}\in M(r,\C)$.
We have $R_{\mathbf{0}} = 0$ where $\mathbf{0} = (0)_i\in \widetilde{\Lambda}$ since $R\in M(r,(\mathfrak{n}\widehat{\mathcal{E}}(\mathfrak{n}))_\Lambda)$.
We have to show the existence of $L$ and $T$.
Write $L = 1_r + \sum_{\mathbf{n}\in\widetilde{\Lambda}}L_{\mathbf{n}}E^{\mathbf{n}}$ and $T = \sum_{\mathbf{n}\in\widetilde{\Lambda}}T_{\mathbf{n}}E^{\mathbf{n}}$.
Then $(H1_r - Q_0 - T)L = L(H1_r - Q_0 - R)$ is equivalent to
\[
	f(\mathbf{n})(H)L_{\mathbf{n}} - [Q_0,L_{\mathbf{n}}] = T_{\mathbf{n}} + S_{\mathbf{n}} - R_{\mathbf{n}},
\]
where $S_{\mathbf{n}}$ is defined by
\[
	\sum_{\mathbf{n}\in\widetilde{\Lambda}}S_{\mathbf{n}}E^{\mathbf{n}} = T(L - 1_r) - (L - 1_r)R. 
\]
By Proposition~\ref{prop:fundamental property of ^(lambda)} the above equation is equivalent to
\[
	\sum_{f(\mathbf{n}) = \mu} S_{\mathbf{n}}E^{\mathbf{n}} = 
	\sum_{f(\mathbf{k}) + f(\mathbf{l}) = \mu}T_{\mathbf{k}}L_{\mathbf{l}}E^{\mathbf{k}}E^{\mathbf{l}} - 
	\sum_{f(\mathbf{k}) + f(\mathbf{l}) = \mu}L_{\mathbf{k}}R_{\mathbf{l}}E^{\mathbf{k}}E^{\mathbf{l}}
\]
for all $\mu\in\mathfrak{a}^*$.
Notice that $L_{\mathbf{0}} = T_{\mathbf{0}} = 0$.
$S_{\mathbf{n}}$ can be defined from the data $\{T_{\mathbf{k}}\mid \mathbf{k} < \mathbf{n}\}$, $\{L_{\mathbf{k}}\mid \mathbf{k} < \mathbf{n}\}$ and $\{R_{\mathbf{k}}\mid \mathbf{k} < \mathbf{n}\}$.

Now we prove the existence of $L$ and $T$.
We choose the $L_{\mathbf{n}}$ and $T_{\mathbf{n}}$ which satisfy
\[
	\begin{cases}
	L_{\mathbf{0}} = 0,\quad T_{\mathbf{0}} = 0,\\
	f(\mathbf{n})(H)L_{\mathbf{n}} - [Q_0,L_{\mathbf{n}}] = T_{\mathbf{n}} + S_{\mathbf{n}} - R_{\mathbf{n}},\\
	\text{if $(Q_0)_{ii} - (Q_0)_{jj}\ne f(\mathbf{n})(H)$ then $(T_{\mathbf{n}})_{ij} = 0$}.
	\end{cases}
\]
By Lemma~\ref{lem:Lemma of linear algebra}, we can choose such $L_{\mathbf{n}}$ and $T_{\mathbf{n}}$ inductively.
Put $L = 1_r + \sum_{\mathbf{n}\in\widetilde{\Lambda}}L_{\mathbf{n}}E^{\mathbf{n}}$ and $T = \sum_{\mathbf{n}\in\widetilde{\Lambda}}T_{\mathbf{n}}E^{\mathbf{n}}$.
Since
\[
	[H,T_{ij}] = \sum_{\mathbf{n}\in\widetilde{\Lambda}}(f(\mathbf{n})(H))(T_{\mathbf{n}})_{ij}E^{\mathbf{n}} = ((Q_0)_{ii} - (Q_0)_{jj})T_{ij},
\]
$L$ and $T$ satisfy the conditions of the lemma.
\end{proof}

\begin{proof}[Proof of Theorem~\ref{thm:definition of boundary value map}]
We can choose positive integers $C = (C_i)\in\Z_{>0}^l$ such that
\[
	\textstyle\{\alpha\in\Lambda^{++}\mid \alpha(\sum_i C_iH_i) = (\lambda_i - \lambda_j)(\sum_i C_iH_i)\} = 
	\begin{cases}
	\{\lambda_i - \lambda_j\} & (\lambda_i - \lambda_j \in \Lambda^{++}),\\
	\emptyset & (\lambda_i - \lambda_j \not\in \Lambda^{++}).
	\end{cases}
\]
The existence of such $C$ is shown by Oshima~\cite[Lemma~2.3]{MR810637}.
Set $X = \sum_i C_iH_i$.
Notice that $(\lambda_i - \lambda_j)(X) \in \{\mu(X)\mid \mu\in\Lambda^{++}\}$ if and only if $\lambda_i - \lambda_j\in \Lambda^{++}$.
By Lemma~\ref{lem:Main lemma of def of boundary value map}, there exist $T\in M(r,(\mathfrak{n}U(\mathfrak{n}))_\Lambda)$ and $L\in M(r,\widehat{\mathcal{E}}(\mathfrak{n})_\Lambda)$ such that
\[
	\begin{cases}
	L \equiv 1_r \pmod{(\mathfrak{n}\widehat{\mathcal{E}}(\mathfrak{n}))_{\Lambda}},\\
	(X1_r - \overline{Q}(X) - T)L = L(X1_r - \overline{Q}(X) - R),\\
	\text{if $\lambda_i - \lambda_j \not\in \Lambda^{++}$ then $T_{ij} = 0$},\\
	\text{if $\lambda_i - \lambda_j \in \Lambda^{++}$ then $[X,T_{ij}] = (\lambda_i - \lambda_j)(X)T_{ij}$}.
	\end{cases}
\]
Let $S\in M(N,(\mathfrak{n}U(\mathfrak{n}))_\Lambda)$ such that $u = \overline{A}w + Su$.
Put $A = (1 - S)^{-1}\overline{A}L^{-1}$, $B = L\overline{B}$ and $v = (v_1,v_2,\dots,v_r) = Bu = Lw$ then $ABu = (1 - S)^{-1}\overline{A}L^{-1}L\overline{B}u = (1 - S)^{-1}\overline{A}w = u$.
Moreover, we have $(X1_r - \overline{Q}(X) - T)v = 0$.
Since $[X,T_{ij}] = (\lambda_i - \lambda_j)(X)T_{ij}$, we have $(X - \lambda_i(X))^rv_i = 0$.

Fix a positive integer $k$ such that $1\le k\le l$.
We can choose a matrix $R_k\in M(r,(\mathfrak{n}\widehat{\mathcal{E}}(\mathfrak{n}))_\Lambda)$ such that $H_kw = (\overline{Q}(H_k) + R_k)w$ by Lemma~\ref{lem:lemma for proof of boudary value map}.
Set $T_k = H_k1_r - \overline{Q}(H_k) - L(H_k1_r - \overline{Q}(H_k) - R_k)L^{-1}$.
Then we have $(H_k1_r - \overline{Q}(H_k) - T_k)v = 0$, i.e.,
\[
	H_kv_i - \sum_{j = 1}^r\overline{Q}(H_k)_{ij}v_j - \sum_{j = 1}^r(T_k)_{ij}v_j = 0
\]
for each $i = 1,2,\dots,r$.
By Corollary~\ref{cor:induced equation}, we have
\[
	H_kv_i - \sum_{j = 1}^r\overline{Q}(H_k)_{ij}v_j - \sum_{j = 1}^r (T_k)_{ij}^{(X,(\lambda_i - \lambda_j)(X))}v_j = 0.
\]
Define $T_k'\in M(r,(\mathfrak{n}U(\mathfrak{n}))_\Lambda)$ by $(T_k')_{ij} = (T_k)_{ij}^{(X,(\lambda_i - \lambda_j)(X))}$.
Since $(T_k)_{ij}^{(\mu)} = 0$ for $\mu\not\in\Lambda^{++}$, we have 
\[
	(T_k)_{ij}^{(X,(\lambda_i - \lambda_j)(X))}
	= \sum_{\mu \in \Lambda^{++},\ \mu(X) = (\lambda_i - \lambda_j)(X)}(T_k)_{ij}^{(\mu)} = 
	\begin{cases}
	(T_k)_{ij}^{(\lambda_i - \lambda_j)} & (\lambda_i - \lambda_j\in\Lambda^{++}),\\
	0 & (\lambda_i - \lambda_j\not\in\Lambda^{++})
	\end{cases}
\]
by the condition of $C$.
In particular $[H,(T_k')_{ij}] = (\lambda_i - \lambda_j)(H)$ for all $H\in \mathfrak{a}$.
Put $Q(\sum x_iH_i) = \overline{Q}(\sum x_iH_i) + \sum x_iT_i'$ for $(x_1,x_2,\dots,x_l)\in \C^l$ then $Q$ satisfies the condition of the theorem.
\end{proof}

Set $\rho = \sum_{\alpha\in\Sigma^+}(\dim \mathfrak{g}_\alpha/2)\alpha$.
From the Iwasawa decomposition $\mathfrak{g} = \mathfrak{k}\oplus\mathfrak{a}\oplus\mathfrak{n}$ we have the decomposition into the direct sum
\[
	U(\mathfrak{g}) = \mathfrak{n}U(\mathfrak{a}\oplus\mathfrak{n})\oplus U(\mathfrak{a})\oplus U(\mathfrak{g})\mathfrak{k}.
\]
Let $\chi_1$ be the projection of $U(\mathfrak{g})$ to $U(\mathfrak{a})$ with respect to this decomposition and $\chi_2$ the algebra automorphism of $U(\mathfrak{a})$ defined by $\chi_2(H) = H - \rho(H)$ for $H\in\mathfrak{a}$.
We define $\chi\colon U(\mathfrak{g})^\mathfrak{k}\to U(\mathfrak{a})$ by $\chi = \chi_2\circ\chi_1$ where $U(\mathfrak{g})^\mathfrak{k} = \{u\in U(\mathfrak{g})\mid \text{$Xu = uX$ for all $X\in \mathfrak{k}$}\}$.
It is known that an image of $U(\mathfrak{g})^\mathfrak{k}$ under $\chi$ is contained in the set of $W$-invariant elements in $U(\mathfrak{a})$.

Fix $\lambda\in \mathfrak{a}^*$.
We can define the algebra homomorphism $U(\mathfrak{a})\to \C$ by $H\mapsto \lambda(H)$ for $H\in \mathfrak{a}$.
We denote this map by the same letter $\lambda$.
Put $\chi_\lambda = \lambda\circ\chi$.
Set $U(\lambda) = U(\mathfrak{g})/(U(\mathfrak{g})\Ker\chi_\lambda + U(\mathfrak{g})\mathfrak{k})$, $u_\lambda = 1\mod{(U(\mathfrak{g})\Ker\chi_\lambda + U(\mathfrak{g})\mathfrak{k})}\in U(\lambda)$ and $U(\lambda)_0 = U(\mathfrak{g})_{2\mathcal{P}}u_\lambda$.
Before applying Theorem~\ref{thm:definition of boundary value map} to $U(\lambda)_0$, we give some lemmas.

\begin{lem}\label{lem:deformation by k}
Let $u\in U(\mathfrak{g})_\mu$ where $\mu \in \mathfrak{a}^*$.
Then there exists an element $x\in U(\mathfrak{g})\mathfrak{k}$ such that $u + x\in U(\mathfrak{a}\oplus\mathfrak{n})_{\mu + 2\mathcal{P}}$.
\end{lem}
\begin{proof}
Set $\overline{\mathfrak{n}} = \theta(\mathfrak{n})$.
We may assume $u\in U(\overline{\mathfrak{n}})_\mu$.
Let $\{U_n(\overline{\mathfrak{n}})\}_{n\in\Z_{\ge 0}}$ be the standard filtration of $U(\mathfrak{n})$ and $n$ the smallest integer such that $u\in U_n(\overline{\mathfrak{n}})$.
We will prove the existence of $x$ by the induction on $n$.

If $n = 0$ then the lemma is obvious.
Assume $n > 0$.
We may assume that there exist a restricted root $\alpha\in\Sigma^+$, an element $u_0\in U_{n - 1}(\overline{\mathfrak{n}})_{\mu + \alpha}$ and a vector $E_{-\alpha}\in \mathfrak{g}_{-\alpha}$ such that $u = u_0 E_{-\alpha}$.
Set $E_\alpha = \theta(E_{-\alpha})$, $u_1 = u_0E_\alpha$, $u_2 = E_\alpha u_0$ and $u_3 = u_1 - u_2$.
Then $u + u_2 + u_3 = u + u_1\in U(\mathfrak{g})\mathfrak{k}$, $u_1,u_2\in U(\mathfrak{g})_{\mu + 2\alpha}$ and $u_3\in U_{n - 1}(\mathfrak{g})_{\mu + 2\alpha}$.
Using the Poincar\'e-Birkhoff-Witt theorem and the induction hypothesis we can choose an element $u_5\in U(\mathfrak{a}\oplus\mathfrak{n})_{\mu + 2\mathcal{P}}$ such that $u_3 - u_5\in U(\mathfrak{g})\mathfrak{k}$.
Again by the induction hypothesis we can choose an element $u_6\in U(\mathfrak{a}\oplus\mathfrak{n})_{\mu + \alpha + 2\mathcal{P}}$ such that $u_0 - u_6\in U(\mathfrak{g})\mathfrak{k}$.
Then $u + u_5 + E_\alpha u_6\in U(\mathfrak{g})\mathfrak{k}$, $u_5\in U(\mathfrak{a}\oplus\mathfrak{n})_{\mu + 2\mathcal{P}}$ and $E_\alpha u_6\in U(\mathfrak{a}\oplus\mathfrak{n})_{\mu + 2\alpha + 2\mathcal{P}}$.
\end{proof}

\begin{lem}\label{lem:lemma for only t^2}
The following equations hold.
\begin{enumerate}
\item $\Ker \chi_\lambda\subset U(\mathfrak{g})_{2\mathcal{P}}$.
\item $U(\mathfrak{a}\oplus\mathfrak{n})\cap (\Ker\chi_\lambda + U(\mathfrak{g})\mathfrak{k}) \subset U(\mathfrak{a}\oplus \mathfrak{n})_{2\mathcal{P}} \cap (\Ker\chi_\lambda + U(\mathfrak{g})\mathfrak{k})$.
\item $U(\mathfrak{a}\oplus\mathfrak{n})\cap (U(\mathfrak{a}\oplus\mathfrak{n})\Ker\chi_\lambda + U(\mathfrak{g})\mathfrak{k})
	\subset U(\mathfrak{a}\oplus \mathfrak{n})(U(\mathfrak{a}\oplus\mathfrak{n})\cap (\Ker\chi_\lambda + U(\mathfrak{g})\mathfrak{k}))$.
\item $U(\mathfrak{a}\oplus \mathfrak{n})\cap (U(\mathfrak{g})\Ker\chi_\lambda + U(\mathfrak{g})\mathfrak{k}) = U(\mathfrak{a}\oplus\mathfrak{n})((U(\mathfrak{g})\Ker\chi_\lambda + U(\mathfrak{g})\mathfrak{k})\cap U(\mathfrak{a}\oplus\mathfrak{n})_{2\mathcal{P}})$.
\end{enumerate}
\end{lem}
\begin{proof}
(1)
It is sufficient to prove $U(\mathfrak{g})^\mathfrak{k} \subset U(\mathfrak{g})_{2\mathcal{P}}$.
Let $G$ be a connected Lie group whose Lie algebra is $\mathfrak{g}_0$ and $K$  its maximal compact subgroup such that $\Lie(K) = \mathfrak{k}_0$.
Since $K$ is connected, $U(\mathfrak{g})^\mathfrak{k} = U(\mathfrak{g})^K = \{u\in U(\mathfrak{g})\mid \text{$\Ad(k)u = u$ for all $k\in K$}\}$.
Assume that $G$ has the complexification $G_\C$.
Fix a maximal abelian subspace $\mathfrak{t}_0$ of $\mathfrak{m}_0$.
Let $K_{\mathrm{split}}$ and $A_{\mathrm{split}}$ be the analytic subgroups with Lie algebras given as the intersections of $\mathfrak{k}_0$ and $\mathfrak{a}_0$ with $[Z_{\mathfrak{g}_0}(\mathfrak{t}_0),Z_{\mathfrak{g}_0}(\mathfrak{t}_0)]$ where $Z_{\mathfrak{g}_0}(\mathfrak{t}_0)$ is the centralizer of $\mathfrak{t}_0$ in $\mathfrak{g}_0$.
Let $F$ be the centralizer of $A_{\mathrm{split}}$ in $K_{\mathrm{split}}$.
Since $F\subset K$, we have $U(\mathfrak{g})^K \subset U(\mathfrak{g})^F$.
On the other hand, we have $U(\mathfrak{g})^F\subset U(\mathfrak{g})_{2\mathcal{P}}$ (See Knapp~\cite[Thorem~7.55]{MR1920389} and Lepowsky~\cite[Proposition~6.1, Proposition~6.4]{MR0379613}).
Hence (1) follows.

(2)
Let $u\in \Ker\chi_\lambda$ and $x\in U(\mathfrak{g})\mathfrak{k}$ such that $u + x \in U(\mathfrak{a}\oplus\mathfrak{n})$.
We can write $u = \sum_\mu u_\mu$ where $u_\mu \in U(\mathfrak{g})_\mu$.
By (1), we have $u_\mu = 0$ for $\mu\not\in 2\mathcal{P}$.
Let $\mu \in 2\mathcal{P}$.
By Lemma~\ref{lem:deformation by k}, there exists an element $y_\mu\in U(\mathfrak{g})\mathfrak{k}$ such that $u_\mu + y_\mu\in U(\mathfrak{a}\oplus\mathfrak{n})_{\mu + 2\mathcal{P}} = U(\mathfrak{a}\oplus\mathfrak{n})_{2\mathcal{P}}$.
Put $y = \sum_\mu y_\mu$.
Then $u + y \in U(\mathfrak{a}\oplus\mathfrak{n})_{2\mathcal{P}}$.
Since $u + y\in U(\mathfrak{a}\oplus\mathfrak{n})$ and $x,y\in U(\mathfrak{g})\mathfrak{k}$ we have $y = x$ by the Poincar\'e-Birkhoff-Witt theorem.
Hence we have $u + x = u + y\in U(\mathfrak{a}\oplus\mathfrak{n})_{2\mathcal{P}}$.

(3)
Let $\sum_i x_iu_i\in U(\mathfrak{a}\oplus\mathfrak{n})\Ker\chi_\lambda$ where $x_i\in U(\mathfrak{a}\oplus\mathfrak{n})$ and $u_i\in \Ker\chi_\lambda$.
We write $u_i = \sum_j z_j^{(i)}v_j^{(i)}$ where $z_j^{(i)}\in U(\mathfrak{a}\oplus\mathfrak{n})$ and $v_j^{(i)}\in U(\mathfrak{k})$.
Let $y\in U(\mathfrak{g})\mathfrak{k}$ and assume $\sum_i x_iu_i + y\in U(\mathfrak{a}\oplus\mathfrak{n})$.
By the Poincar\'e-Birkhoff-Witt theorem, $\sum_i x_iu_i + y = \sum_{i,j} x_iz_j^{(i)}v_{j,0}^{(i)}$ where $v_{j,0}^{(i)}$ is the constant term of $v_j^{(i)}$.
Hence $\sum_i x_iu_i + y = \sum_i x_i(u_i + \sum_j z_j^{(i)}(v_{j,0}^{(i)} - v_j^{(i)})) \in U(\mathfrak{a}\oplus\mathfrak{n})(U(\mathfrak{a}\oplus\mathfrak{n})\cap (\Ker \chi_\lambda + U(\mathfrak{g})\mathfrak{k}))$.

(4)
Since $\Ker\chi_\lambda \subset U(\mathfrak{g})^\mathfrak{k}$, we have
\[
	U(\mathfrak{g})\Ker\chi_\lambda + U(\mathfrak{g})\mathfrak{k} = U(\mathfrak{a}\oplus\mathfrak{n})(\Ker\chi_\lambda)U(\mathfrak{k}) + U(\mathfrak{g})\mathfrak{k} = U(\mathfrak{a}\oplus\mathfrak{n})\Ker\chi_\lambda + U(\mathfrak{g})\mathfrak{k}.
\]
By (2) and (3), we have
\begin{align*}
U(\mathfrak{a}\oplus\mathfrak{n})\cap (U(\mathfrak{g})\Ker\chi_\lambda + U(\mathfrak{g})\mathfrak{k})
& = U(\mathfrak{a}\oplus\mathfrak{n})\cap (U(\mathfrak{a}\oplus\mathfrak{n})\Ker\chi_\lambda + U(\mathfrak{g})\mathfrak{k})\\
& \subset U(\mathfrak{a}\oplus\mathfrak{n})(U(\mathfrak{a}\oplus\mathfrak{n})_{2\mathcal{P}}\cap (\Ker\chi_\lambda + U(\mathfrak{g})\mathfrak{k}))\\
& \subset U(\mathfrak{a}\oplus\mathfrak{n})((U(\mathfrak{g})\Ker\chi_\lambda + U(\mathfrak{g})\mathfrak{k})\cap U(\mathfrak{a}\oplus\mathfrak{n})_{2\mathcal{P}}).
\end{align*}
This implies (4).
\end{proof}

\begin{lem}\label{lem:only t^2}
We have the following equations.
\begin{enumerate}
\item $U(\lambda)_0 = U(\mathfrak{a}\oplus\mathfrak{n})_{2\mathcal{P}}u_\lambda$.
\item $U(\mathfrak{n})\otimes_{U(\mathfrak{n})_{2\mathcal{P}}}U(\lambda)_0 \simeq U(\lambda)$ under the map $p\otimes u\mapsto pu$.
\end{enumerate}
\end{lem}
\begin{proof}
(1) This is obvious from Lemma~\ref{lem:deformation by k}.

(2)
Let $I = U(\mathfrak{a}\oplus \mathfrak{n})\cap (U(\mathfrak{g})\Ker\chi_\lambda + U(\mathfrak{g})\mathfrak{k})$.
We have $U(\mathfrak{a}\oplus\mathfrak{n})\otimes_{U(\mathfrak{a}\oplus\mathfrak{n})_{2\mathcal{P}}}U = U(\mathfrak{n})\otimes_{U(\mathfrak{n})_{2\mathcal{P}}}U$ for any $U(\mathfrak{a}\oplus\mathfrak{n})_{2\mathcal{P}}$-module $U$ since $U(\mathfrak{a}\oplus\mathfrak{n})_{2\mathcal{P}} = U(\mathfrak{a})\otimes U(\mathfrak{n})_{2\mathcal{P}}$.

By (1), we have $U(\lambda)_0 = U(\mathfrak{a}\oplus\mathfrak{n})_{2\mathcal{P}}/(I\cap U(\mathfrak{a}\oplus\mathfrak{n})_{2\mathcal{P}})$.
Hence
\begin{align*}
U(\mathfrak{n})\otimes_{U(\mathfrak{n})_{2\mathcal{P}}}U(\lambda)_0
& = U(\mathfrak{a}\oplus\mathfrak{n})\otimes_{U(\mathfrak{a}\oplus\mathfrak{n})_{2\mathcal{P}}}U(\lambda)_0\\
& = U(\mathfrak{a}\oplus\mathfrak{n})\otimes_{U(\mathfrak{a}\oplus\mathfrak{n})_{2\mathcal{P}}}(U(\mathfrak{a}\oplus\mathfrak{n})_{2\mathcal{P}}/(I\cap U(\mathfrak{a}\oplus\mathfrak{n})_{2\mathcal{P}}))\\
& = U(\mathfrak{a}\oplus\mathfrak{n})/(U(\mathfrak{a}\oplus\mathfrak{n})(I\cap U(\mathfrak{a}\oplus\mathfrak{n})_{2\mathcal{P}}))\\
& = U(\mathfrak{a}\oplus\mathfrak{n})/I\\
& = U(\lambda)
\end{align*}
by Lemma~\ref{lem:lemma for only t^2} (4).
\end{proof}

\begin{lem}\label{lem:soe finiteness property of U(n)_2P}
Let $\{U_n(\mathfrak{n})\}_{n\in\Z_{\ge 0}}$ be the standard filtration of $U(\mathfrak{n})$ and $U_n(\mathfrak{n})_{2\mathcal{P}} = U_n(\mathfrak{n})\cap U(\mathfrak{n})_{2\mathcal{P}}$.
Set $U_{-1}(\mathfrak{n}) = U_{-1}(\mathfrak{n})_{2\mathcal{P}} = 0$, $R = \Gr U(\mathfrak{n})_{2\mathcal{P}} = \bigoplus_n U_n(\mathfrak{n})_{2\mathcal{P}}/U_{n - 1}(\mathfrak{n})_{2\mathcal{P}}$ and $R' = \Gr U(\mathfrak{n}) = \bigoplus_n U_n(\mathfrak{n})/U_{n - 1}(\mathfrak{n})$.
\begin{enumerate}
\item $R'$ is a finitely generated $R$-module.
\item $U(\mathfrak{n})$ is a finitely generated $U(\mathfrak{n})_{2\mathcal{P}}$-module.
\item $U(\mathfrak{n})_{2\mathcal{P}}$ is right and left Noetherian.
\item $U(\lambda)_0$ is a finitely generated $U(\mathfrak{n})_{2\mathcal{P}}$-module.
\end{enumerate}
\end{lem}
\begin{proof}
(1)
Let $\Gamma = \{E^\varepsilon\mid \varepsilon\in\{0,1\}^m\}$.
We denote the principal symbol of $u\in U(\mathfrak{n})$ by $\sigma(u)$.
Notice that if $u\in U(\mathfrak{n})_{2\mathcal{P}}$ then $\sigma(u)$ is the principal symbol of $u$ as an element of $U(\mathfrak{n})_{2\mathcal{P}}$.

We will prove that $\{\sigma(E)\mid E\in \Gamma\}$ generates $R'$ as an $R$-module.
Let $x\in R'$.
We may assume that $x$ is homogeneous, thus there exists an element $u\in U(\mathfrak{n})$ such that $x = \sigma(u)$.
Moreover we may assume that there exist non-negative integers $p = (p_1,p_2,\dots,p_m)$ such that $u = E^p$.
Choose $\varepsilon_i\in \{0,1\}$ such that $\varepsilon_i\equiv p_i\pmod{2}$.
Set $q_i = (p_i - \varepsilon_i)/2 \in \Z_{\ge 0}$, $\varepsilon = (\varepsilon_1,\varepsilon_2,\dots,\varepsilon_m)$ and $q = (q_1,q_2,\dots,q_m)$.
Then we have $x = \sigma(E^p) = \sigma(E^{2q})\sigma(E^\varepsilon)$.
Since $\sigma(E^{2q})\in R$, this implies that $\{\sigma(E)\mid E\in \Gamma\}$ generates $R'$ as an $R$-module.

(2) This is a direct consequence of (1).

(3)
By the Poincar\'e-Birkhoff-Witt theorem, $R'$ is isomorphic to a polynomial ring.
In particular $R'$ is Noetherian.
By the theorem of Eakin-Nagata and (1), we have $R$ is Noetherian.
This implies (3).

(4)
Since $U(\lambda)$ is a finitely generated $U(\mathfrak{n})$-module and $U(\mathfrak{n})$ is a finitely generated $U(\mathfrak{n})_{2\mathcal{P}}$-module, $U(\lambda)$ is a finitely generated $U(\mathfrak{n})_{2\mathcal{P}}$-module.
Hence $U(\lambda)_0$ is a finitely generated $U(\mathfrak{n})_{2\mathcal{P}}$-module by (3).
\end{proof}

We enumerate $W = \{w_1,w_2,\dots,w_r\}$ such that $\re w_1\lambda\ge \re w_2\lambda\ge \dots \ge \re w_r\lambda$.
\begin{thm}\label{thm:definition of boundary value map for G/K}
There exist matrices $A\in M(1,r,\widehat{\mathcal{E}}(\mathfrak{n})_{2\mathcal{P}})$ and $B\in M(r,1,\widehat{\mathcal{E}}(\mathfrak{a}\oplus\mathfrak{n},\mathfrak{n})_{2\mathcal{P}})$ such that $v_\lambda = Bu_\lambda\in(\widehat{\mathcal{E}}(\mathfrak{g},\mathfrak{n})\otimes_{U(\mathfrak{g})}U(\lambda))^r$ satisfies the following conditions:
\begin{itemize}
\item There exists a linear map $Q\colon \mathfrak{a}\to M(r,U(\mathfrak{n})_{2\mathcal{P}})$ such that
\[
	\begin{cases}
	\text{$Hv_\lambda = Q(H)v_\lambda$ for all $H\in\mathfrak{a}$},\\
	\text{$Q(H)_{ii} = (\rho + w_i\lambda)(H)$ for all $H\in\mathfrak{a}$},\\
	\text{if $w_i\lambda - w_j\lambda \not\in2\mathcal{P}^+$ then $Q(H)_{ij} = 0$ for all $H\in\mathfrak{a}$},\\
	\text{if $w_i\lambda - w_j\lambda \in 2\mathcal{P}^+$ then $[H',Q(H)_{ij}] = (w_i\lambda - w_j\lambda)(H')Q(H)_{ij}$ for all $H,H'\in \mathfrak{a}$}.
	\end{cases}
\]
\item We have $u_\lambda = Av_\lambda$.
\item Let $(v_1,v_2,\dots,v_r) = v_\lambda$. Then $\{v_i\pmod{\mathfrak{n}U(\lambda)}\}$ is a basis of $U(\lambda)/\mathfrak{n}U(\lambda)$.
\end{itemize}
\end{thm}
\begin{proof}
Let $u_1,u_2,\dots,u_N$ be generators of $U(\lambda)_0$ as a $U(\mathfrak{n})_{2\mathcal{P}}$-module.
These are also generators of $U(\lambda)$ as a $U(\mathfrak{n})$-module by Lemma~\ref{lem:only t^2}.
We choose matrices $C = {}^t(C_1,C_2,\dots,C_N)\in M(N,1,U(\mathfrak{a}\oplus \mathfrak{n})_{2\mathcal{P}})$ and $D = (D_1,D_2,\dots,D_N)\in M(1,N,U(\mathfrak{n})_{2\mathcal{P}})$ such that ${}^t(u_1,u_2,\dots,u_N) = Cu_\lambda$ and $u_\lambda = D\,{}^t(u_1,u_2,\dots,u_N)$.

Notice that $U(\mathfrak{n})_{2\mathcal{P}} + \mathfrak{n}U(\mathfrak{n}) = U(\mathfrak{n})$.
By Lemma~\ref{lem:only t^2},
\begin{align*}
U(\lambda)/\mathfrak{n}U(\lambda) & = 
(U(\mathfrak{n})/\mathfrak{n}U(\mathfrak{n}))\otimes_{U(\mathfrak{n})}U(\lambda)\\
& = (U(\mathfrak{n})/\mathfrak{n}U(\mathfrak{n}))\otimes_{U(\mathfrak{n})}U(\mathfrak{n})\otimes_{U(\mathfrak{n})_{2\mathcal{P}}}U(\lambda)_0\\
& = (U(\mathfrak{n})/\mathfrak{n}U(\mathfrak{n}))\otimes_{U(\mathfrak{n})_{2\mathcal{P}}}U(\lambda)_0\\
& = ((U(\mathfrak{n})_{2\mathcal{P}} + \mathfrak{n}U(\mathfrak{n}))/\mathfrak{n}U(\mathfrak{n}))\otimes_{U(\mathfrak{n})_{2\mathcal{P}}}U(\lambda)_0\\
& = (U(\mathfrak{n})_{2\mathcal{P}}/(\mathfrak{n}U(\mathfrak{n})\cap U(\mathfrak{n})_{2\mathcal{P}}))\otimes_{U(\mathfrak{n})_{2\mathcal{P}}}U(\lambda)_0\\
& = (U(\mathfrak{n})_{2\mathcal{P}}/(\mathfrak{n}U(\mathfrak{n}))_{2\mathcal{P}})\otimes_{U(\mathfrak{n})_{2\mathcal{P}}}U(\lambda)_0\\
& = U(\lambda)_0/(\mathfrak{n}U(\mathfrak{n}))_{2\mathcal{P}}U(\lambda)_0.
\end{align*}
On the other hand,
\begin{align*}
U(\lambda)/\mathfrak{n}U(\lambda) & = U(\mathfrak{g})/(\mathfrak{n}U(\mathfrak{g}) + U(\mathfrak{g})\Ker\chi_\lambda + U(\mathfrak{g})\mathfrak{k})\\
& = (\mathfrak{n}U(\mathfrak{g}) + U(\mathfrak{a}) + U(\mathfrak{g})\mathfrak{k})/(\mathfrak{n}U(\mathfrak{g}) + U(\mathfrak{g})\Ker\chi_\lambda + U(\mathfrak{g})\mathfrak{k})\\
& = U(\mathfrak{a})/((\mathfrak{n}U(\mathfrak{g}) + U(\mathfrak{g})\Ker\chi_\lambda + U(\mathfrak{g})\mathfrak{k})\cap U(\mathfrak{a})).
\end{align*}
By the definition of $\chi_\lambda$, we have
\[
	(\mathfrak{n}U(\mathfrak{g}) + U(\mathfrak{g})\Ker\chi_\lambda + U(\mathfrak{g})\mathfrak{k})\cap U(\mathfrak{a}) = \sum_{p\in U(\mathfrak{a})^W}U(\mathfrak{a})(\chi_2^{-1}(p) - \lambda(p))
\]
where $U(\mathfrak{a})^W$ is a $\C$-algebra of $W$-invariant elements of $U(\mathfrak{a})$.
By the result of Oshima~\cite[Proposition~2.8]{MR1039854}, the set of eigenvalues of $H\in\mathfrak{a}$ on $U(\mathfrak{a})/(\sum_{p\in U(\mathfrak{a})^W}U(\mathfrak{a})(\chi_2^{-1}(p) - \lambda(p)))$ is $\{(\rho + w\lambda)(H) \mid w\in W\}$ with multiplicities.

We take matrices $A'\in M(N,r,\widehat{\mathcal{E}}(\mathfrak{n})_{2\mathcal{P}})$ and $B'\in M(r,N,\widehat{\mathcal{E}}(\mathfrak{n})_{2\mathcal{P}})$ such that the conditions of Theorem~\ref{thm:definition of boundary value map} hold.
Put $A = DA'$, $B = B'C$ then $A,B$ satisfy the conditions of the theorem.
\end{proof}

\section{Structure of Jacquet modules (regular case)}\label{sec:Structure of Jacquet modules (regular case)}
In this section we assume that $\lambda$ is regular, i.e., $w\lambda \ne \lambda$ for all $w\in W\setminus\{ e\}$.
Let $r = \# W$ and $v_\lambda = (v_1,v_2,\dots,v_r)\in(\widehat{\mathcal{E}}(\mathfrak{g},\mathfrak{n})\otimes_{U(\mathfrak{g})}U(\lambda))^r$ as in Theorem~\ref{thm:definition of boundary value map for G/K}.
Set $\mathcal{W}(i) = \{j\mid w_i\lambda - w_j\lambda\in 2\mathcal{P}^+\}$ for each $i = 1,2,\dots,r$.
\begin{thm}\label{thm:relations of v_i}
We have $Xv_i\in \sum_{j\in\mathcal{W}(i)}U(\mathfrak{g})v_j$ for all $X\in\theta(\mathfrak{n})\oplus\mathfrak{m}$.
\end{thm}

Let $A = {}^t(A^{(1)},A^{(2)},\dots,A^{(r)})$ be as in Theorem~\ref{thm:definition of boundary value map for G/K} and $\overline{A} = {}^t(\overline{A^{(1)}},\overline{A^{(2)}},\dots,\overline{A^{(r)}})$ an element of $M(r,1,\C)$ such that $A^{(i)} - \overline{A^{(i)}} \in \mathfrak{n}\widehat{\mathcal{E}}(\mathfrak{n})$.

\begin{lem}
We have $\overline{A^{(i)}} \ne 0$ for each $i = 1,2,\dots,r$.
\end{lem}
\begin{proof}
Put $\overline{U(\lambda)} = U(\lambda)/\mathfrak{n}U(\lambda)$, $\overline{u_\lambda} = u_\lambda\pmod{\mathfrak{n}U(\lambda)}$ and $\overline{v_i} = v_i\pmod{\mathfrak{n}U(\lambda)}$.
Let $\overline{B} = (\overline{B^{(1)}},\overline{B^{(2)}},\dots,\overline{B^{(r)}})$ be a matrix in $M(1,r,U(\mathfrak{a}))$ such that $\overline{v_j} = \overline{B^{(j)}}\overline{u_\lambda}$.
Then we have $\overline{v_j} = \sum_i \overline{A^{(i)}}\,\overline{B^{(j)}}\overline{v_i}$.
By the regularity of $\lambda$, we have $H\overline{v_j} = (w_j\lambda)(H)\overline{v_j}$ and $H\overline{B^{(j)}}\overline{v_i} = (w_i\lambda)(H)\overline{B^{(j)}}\overline{v_i}$ for all $H\in \mathfrak{a}$.
This implies $\overline{A^{(j)}} \ne 0$ since $\lambda$ is regular.
\end{proof}

\begin{proof}[Proof of Theorem~\ref{thm:relations of v_i}]
Put $f(\mathbf{n}) = \sum_i \mathbf{n}_i\beta_i$ for $\mathbf{n} = (\mathbf{n}_i) \in \Z^m$.
Set $\widetilde{\Lambda} = \{\mathbf{n}\in \Z_{\ge 0}^m\mid f(\mathbf{n}) \in 2\mathcal{P}\}$.
We write $A^{(j)} = \sum_{\mathbf{n}\in\widetilde{\Lambda}} A^{(j)}_{\mathbf{n}} E^{\mathbf{n}}$.
Let $\alpha\in\Sigma^+$ and $E_\alpha\in\mathfrak{g}_\alpha$.
Since $\mathfrak{k}u_\lambda = 0$, we have $(\theta(E_\alpha) + E_\alpha)u_\lambda = 0$.
Hence $(\theta(E_\alpha) + E_\alpha)\sum_j\sum_{\mathbf{n}} A^{(j)}_{\mathbf{n}} E^{\mathbf{n}}v_j = 0$.

By applying Corollary~\ref{cor:induced equation} we have
\[
	\sum_{j = 1}^r\left(\sum_{{\mathbf{n}}\in\widetilde{\Lambda}} A^{(j)}_{\mathbf{n}} (\theta(E_\alpha) + E_\alpha)E^{\mathbf{n}}\right)^{(w_i\lambda - w_j\lambda - \alpha )}v_j = 0
\]
for $i = 1,2,\dots,r$.
On one hand if $w_i\lambda - w_j\lambda\not\in 2\mathcal{P}_+$ then 
\[
	\left(\sum_{{\mathbf{n}}\in\widetilde{\Lambda}} A^{(j)}_{\mathbf{n}} (\theta(E_\alpha) + E_\alpha)E^{\mathbf{n}}\right)^{(w_i\lambda - w_j\lambda - \alpha )} = 0.
\]
On the other hand
\[
	\left(\sum_{{\mathbf{n}}\in\widetilde{\Lambda}} A^{(i)}_{\mathbf{n}} (\theta(E_\alpha) + E_\alpha)E^{\mathbf{n}}\right)^{(-\alpha)} = A^{(i)}_{\mathbf{0}}\theta(E_\alpha).
\]
Hence we have
\[
	A_{\mathbf{0}}^{(i)}\theta(E_\alpha)v_i \in \sum_{j\in \mathcal{W}(i)}U(\mathfrak{g})v_j.
\]
Since $A^{(i)}_{\mathbf{0}} = \overline{A^{(i)}} \ne 0$, we have
\[
	\theta(E_{\alpha})v_i \in \sum_{j\in \mathcal{W}(i)}U(\mathfrak{g})v_j.
\]

Next let $X$ be an element of $\mathfrak{m}$.
By Corollary~\ref{cor:induced equation}, we have
\[
	\sum_{j = 1}^r\left(\sum_{\mathbf{n}\in\widetilde{\Lambda}} A^{(j)}_{\mathbf{n}} XE^{\mathbf{n}}\right)^{(w_i\lambda - w_j\lambda)}v_j = 0.
\]
We can prove $Xv_i \in \sum_{j\in \mathcal{W}(i)}U(\mathfrak{g})v_j$ by the same argument.
\end{proof}

Put $V(\lambda) = \sum_i U(\mathfrak{g})v_i \subset \widehat{\mathcal{E}}(\mathfrak{g},\mathfrak{n})\otimes_{U(\mathfrak{g})}U(\lambda)$.
\begin{cor}\label{cor:V = J(U)}
\[
	V(\lambda) = J(U(\lambda)).
\]
\end{cor}
\begin{proof}
By Theorem~\ref{thm:relations of v_i}, $V(\lambda)$ is finitely generated as a $U(\mathfrak{n})$-module.
By applying Proposition~\ref{prop:generalized result of Goodman and Wallach}, we see that the map
$\widehat{\mathcal{E}}(\mathfrak{g},\mathfrak{n})\otimes_{U(\mathfrak{g})}V(\lambda) \to \prod_{\mu\in\mathfrak{a}^*}V(\lambda)_\mu$ is bijective.
Hence
$\widehat{\mathcal{E}}(\mathfrak{g},\mathfrak{n})\otimes_{U(\mathfrak{g})}V(\lambda)\to
\widehat{\mathcal{E}}(\mathfrak{g},\mathfrak{n})\otimes_{U(\mathfrak{g})}U(\lambda)$
is injective by Proposition~\ref{prop:induced equation}.
This map is also surjective since $v_1,v_2,\dots,v_r$ are generators of $\widehat{\mathcal{E}}(\mathfrak{g},\mathfrak{n})\otimes_{U(\mathfrak{g})}U(\lambda)$.

We have 
$\widehat{\mathcal{E}}(\mathfrak{g},\mathfrak{n})\otimes_{U(\mathfrak{g})}V(\lambda) =
\widehat{\mathcal{E}}(\mathfrak{g},\mathfrak{n})\otimes_{U(\mathfrak{g})}U(\lambda)$.
Since $U(\lambda)$ and $V(\lambda)$ are finitely generated as $U(\mathfrak{n})$-modules, we have
\begin{align*}
\widehat{\mathcal{E}}(\mathfrak{g},\mathfrak{n})\otimes_{U(\mathfrak{g})}U(\lambda) = \widehat{J}(U(\lambda)),\\
\widehat{\mathcal{E}}(\mathfrak{g},\mathfrak{n})\otimes_{U(\mathfrak{g})}V(\lambda) = \widehat{J}(V(\lambda)),
\end{align*}
by Proposition~\ref{prop:algebraic property of E(n)}.
Hence we have $J(U(\lambda)) = J(V(\lambda)) = V(\lambda)$ by Corollary~\ref{cor:Jacquet module of trivial case}.
\end{proof}

Recall the definition of a generalized Verma module.
Set $\overline{\mathfrak{p}} = \theta(\mathfrak{p})$, $\overline{\mathfrak{n}} = \theta(\mathfrak{n})$ and $\rho = \sum_{\alpha\in\Sigma^+}(\dim \mathfrak{g}_\alpha/2)\alpha$.

\begin{defn}[Generalized Verma module]\label{defn:generalized Verma module}
Let $\mu\in\mathfrak{a}^*$.
Define the one-dimensional representation $\C_{\rho + \mu}$ of $\overline{\mathfrak{p}}$ by $(X + Y + Z)v = (\rho + \mu)(Y)v$ for $X\in \mathfrak{m}$, $Y\in\mathfrak{a}$, $Z\in\overline{\mathfrak{n}}$, $v\in \C_{\rho + \mu}$.
We define a $U(\mathfrak{g})$-module $M(\mu)$ by
\[
	M(\mu) = U(\mathfrak{g})\otimes_{U(\overline{\mathfrak{p}})}\C_{\rho + \mu}.
\]
This is called a generalized Verma module.
\end{defn}

Set $V_i = \sum_{j \ge i}U(\mathfrak{g})v_j$.
By the universality of tensor products, any $U(\overline{\mathfrak{p}})$-module homomorphism $\C_{\rho + \mu} \to U$ is uniquely extended to the $U(\mathfrak{g})$-module homomorphism $M(\mu)\to U$ for a $U(\mathfrak{g})$-module $U$.
In particular we have the surjective $U(\mathfrak{g})$-module homomorphism $M(w_i\lambda)\to V_i/V_{i + 1}$.
We shall show that $V_i/V_{i + 1}$ is isomorphic to a generalized Verma module using the character theory.

%The bijectivity of this homomorphism is shown by using the character theory.
%We use the character theory.
Let $G$ be a connected Lie group such that $\Lie(G) = \mathfrak{g}_0$, $K$ its maximal compact subgroup with its Lie algebra $\mathfrak{k}_0$, $P$ the parabolic subgroup whose Lie algebra is $\mathfrak{p}_0$ and $P = MAN$ the Langlands decomposition of $P$ where Lie algebra of $M$ (resp.\ $A$, $N$) is $\mathfrak{m}_0$ (resp.\ $\mathfrak{a}_0$, $\mathfrak{n}_0$).

Since $U(w\lambda) = U(\lambda)$ for $w\in W$, we may assume that $\re\lambda$ is dominant, i.e., $\re\lambda(H_i) \le 0$ for each $i = 1,2,\dots,l$.
By the result of Kostant~\cite[Theorem~2.10.3]{MR0399361}, $U(\lambda)$ is isomorphic to the space of $K$-finite vectors of the non-unitary principal series representation $\Ind_{P}^G (1\otimes\lambda)_K$.
% defined by
%\begin{multline*}
%\Ind_P^G(1\otimes\lambda)_K \\
%= \{f\colon G\to \C\mid \text{$f$ is left $K$-finite, $f(gman) = e^{-(\lambda + \rho)(\log a)}f(g)$ for $g\in G$, $m\in M$, $a\in A$, $n\in N$}\}.
%\end{multline*}
The character of this representation is calculated by Harish-Chandra (See Knapp~\cite[Proposition~10.18]{MR1880691}).
Before we state it, we prepare some notations.
Let $H = TA$ be the maximally split Cartan subgroup,
$\mathfrak{h}_0$ its Lie algebra,
$T = H\cap M$,
$\Delta$ the root system of $H$,
$\Delta^+$ the positive system compatible with $\Sigma^+$,
%$\delta = \sum_{\alpha\in\Delta^+}\alpha/2$,
$\Delta_I$ the set of imaginary roots,
$\Delta_I^+ = \Delta^+ \cap \Delta_I$ and
$\xi_\alpha$ the one-dimensional representation of $H$ whose derivation is $\alpha$ for $\alpha\in\mathfrak{h}^*$.
Under these notations, the distribution character $\Theta_G(U(\lambda))$ of $U(\lambda)$ is as follows;
\[
	\Theta_G(U(\lambda))(ta) = \frac{\sum_{w \in W}\xi_{\rho + w\lambda}(a)}{\prod_{\alpha\in\Delta^+\setminus\Delta_I^+}\lvert1 - \xi_\alpha(ta)\rvert}
	\quad (t\in T,\ a\in A).
\]

We will use the Osborne conjecture, which was proved by Hecht and Schmid~\cite[Theorem~3.6]{MR716371}.
To state it, we must define a character of $J(U)$ for a Harish-Chandra module $U$.
Recall that $J(U)$ is an object of the category $\mathcal{O}'_P$, i.e.,
\begin{enumerate}
\item the actions of $M\cap K$ and $\mathfrak{g}$ are compatible,
\item $J(U)$ splits under $\mathfrak{a}$ into a direct sum of generalized weight spaces, each of them being a Harish-Chandra modules for $MA$,
\item $J(U)$ is $U(\overline{\mathfrak{n}})$- and $Z(\mathfrak{g})$-finite
\end{enumerate}
(See Hecht and Schmid~\cite[(34)Lemma]{MR705884}).
For an object $V$ of $\mathcal{O}'_P$, we define the character $\Theta_P(V)$ of $V$ by
\[
	\Theta_P(V) = \sum_{\mu\in\mathfrak{a}^*}\Theta_{MA}(V_\mu),
\]
where $V_\mu$ is a generalized $\mu$-weight space of $V$.
Let $G'$ be the set of regular elements of $G$.
Set
\begin{gather*}
A^{-} = \{a\in A\mid \text{$\alpha(\log a) < 0$ for all $\alpha\in\Sigma^+$}\},\\
(MA)^{-} = \text{interior of $\left\{g\in MA \Biggm| \text{$\prod_{\alpha\in\Delta^+\setminus\Delta_I^+}(1 - \xi_\alpha(ga)) \ge 0$ for all $a\in A^{-}$}\right\}$ in $MA$}.
\end{gather*}
Then the Osborn conjecture says that $\Theta_G(U)$ and $\Theta_P(J(U))$ coincide on $(MA)^{-}\cap G'$ (See Hecht and Schmid~\cite[(42)Lemma]{MR705884}).
It is easy to calculate the character of a generalized Verma module.
We have
\[
	\Theta_P(M(\mu))(ta) = \frac{\xi_{\rho + \mu}(a)}{\prod_{\alpha\in\Delta^+\setminus\Delta^+_I}(1 - \xi_\alpha(ta))} \quad (t\in T, \ a\in A).
\]
Consequently we have
\[
	\Theta_P(J(U(\lambda))) = \sum_{w\in W}\Theta_P(M(w\lambda)).
\]
This implies the following theorem when $\lambda$ is regular.
\begin{thm}\label{thm:Main theorem for regular case}
There exists a filtration $0 = V_{r + 1} \subset V_r \subset \dots \subset V_1 = J(U(\lambda))$ of $J(U(\lambda))$ such that $V_i / V_{i + 1}$ is isomorphic to $M(w_i\lambda)$ for an arbitrary $\lambda\in\mathfrak{a}^*$.
Moreover if $w\lambda - \lambda\not\in 2\mathcal{P}$ for all $w\in W\setminus \{e\}$ then $J(U(\lambda)) \simeq \bigoplus_{w\in W}M(w\lambda)$.
\end{thm}

\section{Structure of Jacquet modules (singular case)}\label{sec:Structure of Jacquet modules (singular case)}
In this section, we shall prove Theorem~\ref{thm:Main theorem for regular case} in the singular case using the translation principle.
%To prove Theorem~\ref{thm:Main theorem for regular case} in the singular case, we use the translation principle.
We retain notations in Section~\ref{sec:Structure of Jacquet modules (regular case)}.
Let $\lambda'$ be an element of $\mathfrak{a}^*$ such that following conditions hold:
\begin{itemize}
\item The weight $\lambda'$ is regular.
\item The weight $(\lambda - \lambda')/2$ is integral.
\item The real part of $\lambda'$ belongs to the same Weyl chamber which real part of $\lambda$ belongs to.
\end{itemize}
First we define the translation functor $T_{\lambda'}^\lambda$.
Let $U$ be a $U(\mathfrak{g})$-module which has an infinitesimal character $\lambda'$. (We regard $\mathfrak{a}^*\subset \mathfrak{h}^*$.)
We define $T_{\lambda'}^\lambda(U)$ by $T_{\lambda'}^\lambda(U) = P_\lambda(U\otimes E_{\lambda - \lambda'})$ where:
\begin{itemize}
\item $E_{\lambda - \lambda'}$ is the finite-dimensional irreducible representation of $\mathfrak{g}$ with an extreme weight $\lambda - \lambda'$.
\item $P_\lambda(V) = \{v\in V\mid \text{for some $n > 0$ and all $z\in Z(\mathfrak{g})$, $(z - \lambda(\widetilde{\chi}(z)))^nv = 0$}\}$ where $Z(\mathfrak{g})$ is the center of $U(\mathfrak{g})$ and $\widetilde{\chi}\colon Z(\mathfrak{g})\to U(\mathfrak{h})$ is the Harish-Chandra homomorphism.
\end{itemize}

Notice that $P_\lambda$ and $T_{\lambda'}^\lambda$ are exact functors.
Theorem~\ref{thm:Main theorem for regular case} in the singular case follows from following two equations.
\begin{enumerate}
\item $T_{\lambda'}^\lambda(U(\lambda')) = U(\lambda)$.
\item $T_{\lambda'}^\lambda(M(w\lambda')) = M(w\lambda)$.
\end{enumerate}

The following lemma is important to prove these equations.
\begin{lem}\label{lem:lemma for weight (by Vogan)}
Let $\nu$ be a weight of $E_{\lambda - \lambda'}$ and $w\in W$.
Assume $\nu = w\lambda - \lambda'$.
Then $\nu = \lambda - \lambda'$.
\end{lem}
\begin{proof}
See Vogan~\cite[Lemma~7.2.18]{MR632407}.
\end{proof}

\begin{proof}[Proof of\/ $T_{\lambda'}^\lambda(U(\lambda')) = U(\lambda)$]
We may assume that $\lambda'$ is dominant.
Notice that we have $U(\lambda') \simeq \Ind_P^G(1\otimes \lambda')_K$.
Let $0 = E_0 \subset E_1 \subset E_2\subset \dots\subset E_n = E_{\lambda - \lambda'}$ be a $P$-stable filtration with the trivial induced action of $N$ on $E_i/E_{i - 1}$.
We may assume that $E_i / E_{i - 1}$ is irreducible.
Let $\nu_i $ be the highest weight of $E_i / E_{i - 1}$.
Notice that $\Ind_P^G(1\otimes \lambda')\otimes E_{\lambda - \lambda'} = \Ind_P^G((1\otimes \lambda')\otimes E_{\lambda - \lambda'})$.
Then $\Ind_P^G(1\otimes \lambda')\otimes E_{\lambda - \lambda'}$ has a filtration $\{M_i\}$ such that $M_i/M_{i - 1} \simeq \Ind_P^G((1\otimes \lambda')\otimes (E_i/E_{i - 1}))$.
Since $\Ind_P^G((1\otimes \lambda')\otimes (E_i/E_{i - 1}))$ has an infinitesimal character $\lambda + \nu_i$, $P_\lambda(M_i/M_{i - 1}) = 0$ if $\nu_i \ne w\lambda - \lambda'$ for all $w\in W$.
By Lemma~\ref{lem:lemma for weight (by Vogan)} we have $T_{\lambda'}^\lambda(\Ind_P^G(1\otimes \lambda') = \Ind_P^G((1\otimes \lambda')\otimes (E_i/E_{i - 1}))$ where $\nu_i = \lambda - \lambda'$.
By the conditions of $\lambda'$, the action of $M$ on $(\lambda - \lambda')$-weight space of $E_{\lambda - \lambda'}$ is trivial.
Consequently $T_{\lambda'}^\lambda(\Ind_P^G(1\otimes \lambda')) = \Ind_P^G((1\otimes\lambda')\otimes(\lambda - \lambda')) = \Ind_P^G(1\otimes \lambda)$.
\end{proof}

\begin{proof}[Proof of\/ $T_{\lambda'}^\lambda(M(w\lambda')) = M(w\lambda)$.]
We may assume $w = e\in W$.
Since $M(\lambda')\otimes E_{\lambda - \lambda'} = U(\mathfrak{g})\otimes (\C_{\lambda'}\otimes E_{\lambda - \lambda'})$, 
the equation follows by the same argument of the proof of $T_{\lambda'}^\lambda(U(\lambda')) = U(\lambda)$.
\end{proof}

%\bibliographystyle{my_amsalpha}
%\bibliography{../../bunken,../../book}
\def\cprime{$'$} \def\dbar{\leavevmode\hbox to 0pt{\hskip.2ex \accent"16\hss}d}

\end{document}